\documentclass[twoside]{article}
\usepackage{graphicx}
\newlength{\ei}\ei=0.0138888889em   
   \newlength{\breit}
\newlength{\SyW}  \newlength{\msu}  \msu=\mathsurround 
\newcommand{\weg}{\hspace{\msu}}\newcommand{\ran}{\hspace{-\mathsurround}}
\newcommand{\qed}{\nopagebreak\hspace*{\fill}\mbox{\scriptsize\rm
  {/}\hspace{-13\ei}{/}\hspace{-13\ei}{/}\hspace{-13\ei}{/}}}   
\newcommand{\nz}{\normalsize}\newcommand{\sz}{\small}

\newcommand{\Ds}{\displaystyle} \newcommand{\Ts}{\textstyle}
  
\newcommand{\Nsum}{\sum\nolimits}

\newcommand{\Tfrac}[2]{{\Ts\frac{#1}{#2}}}

\newcommand{\wsup}{\mathop{\rm sup\vphantom{f}}\nolimits}
\newcommand{\winf}{\mathop{\rm inf\vphantom{p}}\nolimits}
\newcommand{\ggrto}{\relbar\joinrel\longrightarrow}

\newcommand{\gwto}{\mathrel{\mathsurround0em \mbox{$\ggrto$}%
  \llap{\settowidth{\SyW}{$\ggrto$}
  \raisebox{-.15ex}{\makebox[\SyW]{\scriptsize\rm w}}}}}
\newfont{\tenmsam}{msam10}                                 

\font\hMtenrm=cmr10 scaled \magstephalf 
\font\hMtenbf=cmbx10 scaled \magstephalf 
\newcommand{\maxev}{\mathop{\rm max\hspace*{6\ei}ev}\nolimits}

\newcommand{\tr}{\mathop{\rm tr}\nolimits}
\newcommand{\rk}{\mathop{\rm rk}\nolimits}

\newcommand{\hJc}{\mathop{\mbox{\hMtenbf I}}\nolimits}
\newcommand{\lJc}{\mathop{\mbox{\large\bf I}}\nolimits}
\newcommand{\Lo}{\mathop{\rm {{}o{}}}\nolimits} 
\newcommand{\diag}{\mathop{\rm diag}\nolimits}
\newcommand{\lin}{\mathop{\rm lin}\nolimits}
\newcommand{\cl}{\mathop{\it c\ell}\nolimits}
\newcommand{\clin}{\mathop{\cl\lin}\nolimits}
\newcommand{\cqlin}{\mathop{{\cl\hskip6\ei}'\lin}\nolimits} 

\newcommand{\Ew}{\mathop{\rm {{}E{}}}\nolimits} 
\newcommand{\hEw}{\mathop{\mbox{\hMtenrm E}}\nolimits}  
\newcommand{\Cov}{\mathop{\rm {{}C{}}}\nolimits}    
\newcommand{\MSE}{\mathop{\rm MSE}\nolimits}    
\newcommand{\relMSE}{\mathop{\rm relMSE}\nolimits}    

\makeatletter                  
\def\magstepams#1{\ifcase#1 \magstep0\or\magstephalf\or\magstep1\fi\relax}
\font\tenmsbm=msbm10 scaled \magstepams\@ptsize
\font\sevmsbm=msbm7  scaled \magstepams\@ptsize
\def\relaxnext@{\let\next\relax}   \def\noaccents@{\def\accentfam@{0}}
\newfam\msbmfam \textfont\msbmfam=\tenmsbm \scriptfont\msbmfam=\sevmsbm
\def\Bberr{Use Black board bold only in math mode}
\def\Bb{\relaxnext@\ifmmode\let\next\Bb@\else
        \def\next{\errmessage{\Bberr}}\fi\next}
\def\Bb@#1{{\Bb@@{#1}}}\def\Bb@@#1{\noaccents@\fam\msbmfam#1}\makeatother
\newcommand{\EM}{{\mbox{\tenmsbm I}}}       
\newcommand{\R}{{\mbox{\tenmsbm R}}}        
\newcommand{\B}{{\mbox{\tenmsbm B}}}        
\newtheorem{Thm}{Theorem}[section]   
\newtheorem{Prop}[Thm]{Proposition}  \newtheorem{Lem}[Thm]{Lemma}
\newtheorem{Cor}[Thm]{Corollary}     
\newtheorem{Exa}[Thm]{Example}       \newtheorem{Rem}[Thm]{Remark}
\newenvironment{Bew}{\begin{trivlist}\item[]{\sc Proof}}{\end{trivlist}}

\newcommand{\rfi}[1]{\makebox[\parindent][l]{%
                     \makebox[0em][r]{\rm(}\sf#1\rm)}}
\newcounter{ABCc}\renewcommand{\theABCc}{\alph{ABCc}}
\newenvironment{ABC}{\begin{list}{
  \rfi{\theABCc}}{\usecounter{ABCc} \topsep 0ex \partopsep 0ex \itemsep0ex
  \parsep=\parskip \leftmargin 0em \rightmargin 0em \itemindent=\parindent
  \listparindent=\parindent  \labelsep 0em \labelwidth 0em }}{\end{list}}
\newcounter{ABCca}\renewcommand{\theABCca}{\arabic{ABCca}}
\newenvironment{ABCa}{\begin{list}{
  \rfi{\theABCca}}{\usecounter{ABCca} \topsep 0ex \partopsep 0ex \itemsep0ex
  \parsep=\parskip \leftmargin 0em \rightmargin 0em \itemindent=\parindent
  \listparindent=\parindent  \labelsep 0em \labelwidth 0em }}{\end{list}}

\begin{document}
\title{\LARGE\bf Neighborhoods as Nuisance Parameters?
                 Robustness vs.~Semiparametrics.}
\author{\nz\rm Helmut Rieder\\
         \sz\sl University of Bayreuth, Germany}
\date{\sz\rm 24~September~2000, under revision} 
\maketitle
\begin{abstract}\noindent
Deviations from the center within a robust neighborhood
may naturally be considered an infinite dimensional nuisance parameter.
Thus, the semiparametric method may be tried, which is to compute the
scores function for the main parameter minus its orthogonal projection
on the closed linear tangent space for the nuisance parameter,
and then rescale for Fisher consistency.
We derive such a semiparametric influence curve by nonlinear projection
on the tangent balls arising in robust statistics.
\par
This semiparametric influence curve is then compared with the optimally
robust influence curve that minimizes maximum weighted mean square error
of the corresponding asymptotically linear estimators over
infinitesimal neighborhoods.
While there is coincidence for Hellinger balls, at least clipping is
achieved for total variation and contamination neighborhoods, but the
semiparametric method in general falls short to solve the minimax MSE
estimation problem for the gross error models.
\par
The semiparametric approach is carried further to testing contaminated
hypotheses. In the one-sided case, for testing hypotheses defined by any
two closed convex sets of tangents, a saddle point is furnished by
projection on the set of differences of these sets.
For total variation and contamination neighborhoods, we thus recover
the robust asymptotic tests based on least favorable pairs.
\par \vspace{1\smallskipamount plus 10ex}\noindent
\par \noindent
{\sl AMS/MSC-2000 classification}: 62F35, 62G35.
\end{abstract}
\section{Introduction}                \setcounter{equation}{0}\label{s.I}
Robustness and semiparametrics are of the same origin, namely,
the desire to get rid of a narrow parametric model.
The first achieves this by including certain neighborhoods, the second
by introducing possibly infinite dimensional nuisance parameters.
Despite of similar goals, however, the relations between these two modern
statistical developments have not been investigated systematically. %
\par     Three aspects appear to be of conceptual interest. 
\par     \bigskip \noindent    {\bf Nonrobustness of Adaptive Procedures}
Informal statements by Huber (1981; Section~1.2) and
                             (1996; Sections 19,~28)
and similar remarks by Hampel et al.~(1986) 
indicate nonrobustness of adaptive estimators, that is,
of estimators which are asymptotically efficient for the location problem
with unknown symmetric density.
\par
Robustness or not of adaptive estimators for more general semiparametric
models has also been addressed, and declared a field of future research,
by Bickel et al.~(1993; Introduction, p~4). In this remark, specific
reference is made to the infinitesimal setup of Hampel et al.~(1986).
\par     
Section~\ref{s.O} below presents the readily available argument to prove
the nonrobustness in the infinitesimal setup, and describes further
aspects and ensuing problems.
But the thrust of the paper is not in this direction. 
\par     \medskip \noindent    {\bf Common Local Asymptotic Basis}
Semiparametric theory, as treated in the accounts by
Bickel et al.~(1993) and van der Vaart~(1998),
and infinitesimal robustness in fact share the same mathematical basis:
the local asymptotic statistical theory due to LeCam.
\par
This mathematical framework has been explicitly declared in these
monographs, but hidden from the beginning by Hampel et al.~(1986); already
their Definition~1 (Section~2.1, p~84) of the basic notion---influence
function, as some kind of G\^ateaux--derivative---is kept informal,
deprived of neighborhoods and of estimators and their laws.
However, once these ingredients are accounted for mathematically, as
in Rieder~(1994), the similarity becomes obvious, and concerns 
the main issues: derivation of asymptotic lower bounds for estimator risk,
and construction of optimal estimators achieving the bounds.
\par     
In this paper, 
we only argue with risks and neglect estimator constructions.
\par     \medskip {\bf Model Deviations as Nuisance Parameter}
While adaptive procedures are still open for robustification,
another view opens robustness to semiparametric arguments:
In principle, and quite naturally, model deviations can be viewed an
infinite dimensional nuisance parameter. Thus a neighborhood model about
a parametric family may be interpreted as a semiparametric model.
\par   
And this is the thrust of our paper: Bring the semiparametric methods
to bear on robust neighborhood models; in particular, derive the
`efficient' influence curves.
If everything works out the nicest way, the corresponding estimators ought
to be optimal for the parameter of the center model, in the presence of
model deviations; which is exactly what one expects from optimally robust
estimators. 
\par     \medskip\noindent    {\bf Outline of the Paper}
In Section~\ref{s.P}, we set up the semiparametric approach, including a
smooth semiparametric model, tangent sets, influence curves of functionals
and asymptotically linear estimators, projections on tangent spaces,
the canonical influence curve, the Cram\'er--Rao and more general
asymptotic minimax bounds.
The development derives the semiparametric recipe which amounts to compute
the scores function for the main parameter minus its orthogonal projection %
on the closed linear tangent space for the nuisance parameter, 
and then rescale for Fisher consistency.
Our presentation differs from Bickel et al.~(1993; Chapters 2--3)
and van der Vaart~(1998; Chapter~25) in that we simultaneously consider
the whole set of influence curves, not only the canonical one.
Moreover, adaptivity, existence of bounded influence curves,
the case of a finite dimensional nuisance tangent space,
and an asymptotic confidence bound that uses nonlinear projection
on closed convex cones are addressed.
\par  In Section~\ref{s.B}, we formulate the robust setup of
infinitesimal neighborhoods. The embedding in the semiparametric mold
creates an identifiability problem, and requires the so-called idealistic
attitude towards robustness. The determination of the corresponding tangent
sets then leads to balls which span the entire~$L_2$ of expectation zero; in
particular, the subtraction of the projection would annihilate the scores.
But the semiparametric recipe seems intuitively plausible 
even if the nuisance tangent set happens not to be a linear space.
\par  In Section~\ref{s.C}, therefore, we deviate from the dogmatic recipe
and, from the scores, subtract only its projection on the closed balls
themselves. In the Hellinger case, the resulting semiparametric influence
curves conicides with the classically optimum one. In the total variation
and contamination cases, the semiparametric influence curves turn out
clipped versions of the scores. Thus, essential features of the
optimally robust influence curves for these models seem to be
recovered by our nonlinearly modified semiparametric approach.
\par In Section~\ref{s.R}, the semiparametric influence curve is checked
more quantitatively under a specified risk, namely, by comparison with
the optimally robust influence curve that minimizes maximum weighted
mean square error of asymptotically linear estimators over
shrinking neighborhoods.
For Hellinger balls, the two influence curves coincide
(with the classically optimum one).
In the case of total variation balls, the semiparametric influence curve
solves the robust mean square error problem only for a particular bias weight,
respectively, for bias weight one and a different neighborhood radius
(in an example shown to be larger than the given radius).
For parameter dimension one, the comparison is also done with respect
to a certain confidence risk. In the case of contamination neighborhoods,
the semiparametric influence curve is bounded only from above.
As the semiparametric influence curve interchanges linear combination
and truncation, in comparison with the optimally robust influence curve,
the discrepancy between the two seems to increase
with the parameter dimension.
\par  Section~\ref{s.O} gives the brief argument that adaptive estimators
not only inherit the asymptotic efficiency of the estimator they adapt
but also its nonrobustness against infinitesimal gross errors,
if only the canonical influence curve is unbounded. The problem of
adapting optimally robust estimators (estimating out the unknown
nuisance parameter) is suggests itself but not treated further. %
\par  Section~\ref{s.T} applies the ideas to testing.
Thus, the semiparametric extension of Neyman's $C(\alpha)$-tests
(to the case of an infinite dimensional nuisance parameter) may be
modified nonlinearly to become applicable to the robust tangent balls.
Then, for contaminated one- and multisided hypotheses about the main
parameter, at least sensibly bounded test statistics are obtained.
In general, optimally robust tests are not even available, against
which these semiparametric competitors might be judged.
\par  In the one-sided, one parameter case however, and for total variation
and contamination neighborhoods, the asymptotic tests based on least
favorable pairs in the sense of Huber and Strassen~(1973) define the
ultimate robustness standard.
In Section~\ref{s.S}, our semiparametric recipe is able to recover these
optimally robust tests. More generally, a saddle point is furnished for
testing hypotheses defined by any two closed convex sets of tangents,
via (nonlinear) projection on the set of differences of these sets.
\par     \bigskip\noindent    {\bf Conclusions}
In the infinitesimal robust setup, the modified semiparametric recipe mostly
yields estimators and tests which are reasonably robust. \mbox{Exact}
agreement with an optimally robust procedure may be achieved for special
loss functions. But, since there is a difference in general, not all aspects
of the robust model seem to be caught correctly by the semiparametric method.
\mbox{Coincidence} occurs rather in the context of testing than estimation,
and then rather for Hellinger and total variation balls than
for contamination neighborhoods.
The semiparametric method copes with parameter dimension greater than one
more easily than the robust method. 
\section{The Semiparametric Setup}    \setcounter{equation}{0}\label{s.P}
To set up the standard semiparametric framework, we employ some
family~${\cal Q}$ in the set~${\cal M}$ of all probabilities on
some sample space~$(\Omega, {\cal B})$,
\begin{equation}
  {\cal Q} =   \{\, Q_{\theta,\nu} \mid \theta\in \Theta\weg,\:
  \nu\in H_{\theta}\,\}   \subset {\cal M}
\end{equation}
The parameter~$\theta$ of interest is finite ($k$-)dimensional,
out of some open parameter set $\Theta \subset \nolinebreak \R^k$,
whereas~$\nu$ acts as nuisance parameter.
For each~$\theta$, $\nu$~ranges over some set~$H_{\theta}$; typically,
subsets of some infinite dimensional function spaces; densities or
differences of densities (Section~\ref{s.B}).
The observations are assumed independent identically distributed,
  $x_1,\ldots,x_n\sim Q_{\theta,\nu}$.
Estimators of~$\theta$ may be any functions $S_n \colon \Omega^n\to\R^k$
which are product measurable ${\cal B}^n$/Borel~$\B^k$.
Let us fix~$(\theta_0,\nu_0)$, the true but unknown values of main
and nuisance parameter.
\par     Optimality results for the estimation of~$\theta_0$ can in general
only be derived asymptotically, for sample size~$n\to \infty$.
Moreover, to obtain meaningful results, estimators, now estimator sequences
  $S=(S_n)$, must be judged not unly at~$(\theta_0,\nu_0)$ but under
local alternatives about~$(\theta_0,\nu_0)$.
Subsequently, the fixed parameter will be omitted whenever feasible.
Thus, we put $Q_{\theta_0,\nu_0}=Q$, and denote expectation and covariance
under~$Q$ by~$\Ew$ and~$\Cov$. 
Also the spaces~$L_2$ and~$L_{\infty}$ of square integrable and
essentially bounded real functions, respectively, refer to the
fixed $Q=Q_{\theta_0,\nu_0}$.
The corresponding spaces of $\R^k$~valued functions are denoted by
$L_2^k$ and~$L_{\infty}^k$. 
\par     For the local asymptotics a certain smoothness of the parametric
model is required, in the sense of mean square differentiability
at~$(\theta_0,\nu_0)$ of square root densities:
There exists some function $\Lambda\in L_2^k$---the scores function
for the main parameter~$\theta$ at~$(\theta_0,\nu_0)$---such that
for each $a\in\R^k$, and for each $g\in\partial_2 {\cal Q}$ there
is some path $t \mapsto \nu^g_t\in H_{\theta_0+ta}$ such that,
as $t\to0$ in~$\R$,
\begin{equation} \label{e.L2diff}
  \sqrt{dQ_{\theta_0+ta, \hspace{6\ei}\nu^g_t}}=
  \bigl(1+ \Tfrac{1}{2} \hspace{6\ei}t \hspace{6\ei}
  (a' \Lambda+g)\bigr) \sqrt{dQ_{\theta_0,\nu_0}}\,   + \Lo(t)
\end{equation}
In this context, the tangent set
  $\partial {\cal Q}= \partial_1 {\cal Q}+ \partial_2 {\cal Q}$
of the model~${\cal Q}$ at~$(\theta_0,\nu_0)$ appears, where
  $\partial_1 {\cal Q}=\{\,a' \Lambda\mid a\in\R^k\,\}$ is
the tangent space for the first parameter component, 
and $\partial_2 {\cal Q}\subset L_2$ denotes the tangent set
for the nuisance component;
all tangents in either class~$\partial_* {\cal Q}$ necessarily
have expectation zero. The covariance ${\cal I}=\Cov\Lambda$ is
the Fisher information of the $\nu_0$-section~${\cal Q}_{\nu_0}$
of model~${\cal Q}$ for the parameter~$\theta$ at~$\theta_0$;
   $ \rk {\cal I} = k $, by~(\ref{e.c.J}) and~(\ref{e.p.rkJ=k}) below.
\par     As for complete technical details, maybe in slightly different
notations, the reader may consult the textbooks 
by Bickel et al.~(1993; Chapters 2--3), 
van der Vaart~(1998; Chapter~25),
and Rieder\footnote{HR, henceforth}(1994; Chapters 2--4).
\par     Influence functions, or influence curves (IC),~$\psi$ for
model~${\cal Q}$ at~$(\theta_0,\nu_0)$ are defined by the conditions
\begin{equation} \label{e.IC}
   \psi\in L_2^k,\quad
   \Ew \psi=0\weg,\quad \Ew \psi \Lambda'=\EM_k\weg, \quad
   \Ew \psi g =0\enskip \forall g\in \partial_2 {\cal Q} 
\end{equation}
where $\EM_k $ denotes the $k\times k$ identity matrix.
The set of all influence curves for model~${\cal Q}$ at~$(\theta_0,\nu_0)$
is denoted by~$\Psi=\Psi_{\theta_0,\nu_0}$.
\par     On the one hand, influence curves go with functionals
      $T \colon {\cal Q}\to\R^k$ which are differentiable, with respect to
model~${\cal Q}$ at~$(\theta_0,\nu_0)$ in accordance with~(\ref{e.L2diff}),
and are Fisher consistent for the main parameter such that
\begin{equation} \label{e.Tdiff}
   T(Q_{\theta_0+ta,\nu^g_t}) =
   T(Q_{\theta_0,\nu_0}) + \hEw \psi (a' \Lambda+g)\,t + \Lo(t)
   = \theta_0 + t \hspace{6\ei}a + \Lo(t)
\end{equation}
On the other hand, influence curves go with asymptotically linear
estimators. These are estimators $S=(S_n)$ that have an expansion
\begin{equation} \label{e.Sal}
   \sqrt{n}\,(S_n- \theta_0)= \frac{1}{\sqrt{n}\,}
   \sum_{i=1}^n \psi(x_i) \: + \Lo_{{Q^n}}(n^0)
\end{equation}
where the remainder tends to zero in probability, under the sequence
of product measures~$Q^n$. Such estimators are asymptotically normal
in accordance with~(\ref{e.L2diff}): Setting
  $ Q_n(a,g)=Q_{\theta_0+s_na,\,\nu_{s_n}^{\hspace{6\ei}g}} $ for
  $s_n=1/\!\sqrt{n}\,$, their distributions under~$Q_n^n(a,g)$ converge
weakly as $n\to \infty$, for every $a\in\R^k$ and $g\in \partial_2 {\cal Q}$,
\begin{equation} \label{e.ufoN}
   \sqrt{n}\,(S_n- \theta_0)(Q_n^n(a,g)) \gwto {\cal N}(a, \Cov\psi)
\end{equation}
Given any $\psi\in \Psi$,
  $ T(M)=\theta_0+2\int \psi\,\sqrt{dQ}\,\sqrt{d\vphantom{Q}\smash{M}}\,$
and $ S_n=\theta_0+1/n\sum \psi(x_i) $ are 
constructions to achieve (\ref{e.Tdiff}) and~(\ref{e.Sal}), which however
depend on~$(\theta_0,\nu_0)$.
\par     For either tangent set~$\partial_\diamond {\cal Q}$
let $\lin\partial_\diamond {\cal Q}$ and $\clin\partial_\diamond {\cal Q}$
denote the linear span, respectively the closed linear span,
of~$\partial_\diamond {\cal Q}$ in~$L_2$.
Thus, $\clin \partial_1 {\cal Q}=\partial_1 {\cal Q}$, and
$ \clin \partial {\cal Q}=\partial_1 {\cal Q}+\clin \partial_2 {\cal Q}$
as $\dim\partial_1 {\cal Q}$ is finite.
Introduce the orthogonal projection
    $\pi_\diamond \colon L_2\to \clin\partial_\diamond {\cal Q}$
on~$\clin\partial_\diamond {\cal Q}$,
and $\Pi_\diamond \colon L_2^k\to (\clin\partial_\diamond {\cal Q})^k$
the orthogonal projection in the product space; then
    $\Pi_\diamond =(\pi_\diamond,\ldots,\pi_\diamond)'$,
acting coordinatewise.
\par     In view of~(\ref{e.IC}), the projection~$\Pi(\psi)$
on~$(\clin\partial {\cal Q})^k$ must be the same for every
  $\psi\in \Psi$---the shortest, or canonical, influence curve~$\varrho$.
In fact,
\begingroup \mathsurround0em\arraycolsep0em \begin{eqnarray}
\label{e.cIC}    \Pi(\psi)&{}={}& \varrho  =
   {\cal J}^{-1} \bigl(\Lambda- \Pi_2(\Lambda)\big)
   \qquad \forall\,\psi\in \Psi \hspace{-3em} \\
\noalign{\noindent where\nopagebreak} \label{e.c.J}
    {\cal J} &{}={}& \Cov\bigl(\Lambda- \Pi_2(\Lambda)\big)  =
               {\cal I} - \Cov \Pi_2(\Lambda)
\end{eqnarray}\endgroup
denotes the Fisher information of model~${\cal Q}$ for the
parameter~$\theta$ at~$(\theta_0,\nu_0)$.
\par     A little argument shows that the existence of influence curves
is equivalent to regularity, that is, positive definiteness, of~${\cal J}$,
\begin{equation} \label{e.p.rkJ=k}
\mathsurround0em\arraycolsep0em \begin{array}{r@{{}\iff{}}l}
\Ds \Psi\ne \emptyset & \Ds {\cal J} > 0 \\
\rule{0pt}{2.75ex} & \Ds a' \Lambda \notin \clin \partial_2 {\cal Q}
  \quad \forall\,a\in\R^k\weg,\: a\ne0
\end{array}\end{equation}
which condition we want to assume subsequently.
\begin{Rem}\rm [\,adaptivity\,]\enskip  With the nuisance parameter~$\nu$
fixed to~$\nu_0$, the $\nu_0$-section ${\cal Q}_{\nu_0}$ of model~${\cal Q}$
is a model without nuisance parameter,
\begin{equation} \label{e.Qn0}
   {\cal Q}_{\nu_0} = \{\,Q_{\,\theta,\nu_0}\mid \theta\in \Theta\,\}
\end{equation}
satisfying~(\ref{e.L2diff}) with $\partial_2 {\cal Q}_{\nu_0}= \{0\}$ and
$\partial {\cal Q}_{\nu_0}= \partial_1 {\cal Q}$.
Consequently, the canonical influence curve 
and the Fisher information of model~${\cal Q}_{\nu_0}$ for the
\mbox{parameter}~$\theta$ at~$\theta_0$ are given by, respectively,
\begin{equation} \label{e.rhoP}
  \hat{\varrho} = {\cal I}^{-1}\Lambda \weg, \qquad {\cal I}=\Cov\Lambda
\end{equation}
The following bound of~${\cal J}$ by~${\cal I}$,
in the positive definite sense, is an immediate consequence
of (\ref{e.cIC}), (\ref{e.c.J}), and~(\ref{e.rhoP}),
\begin{equation} \label{e.Cadapt} \hspace{-.5em}
  \Cov\hat\varrho= {\cal I}^{-1}\le {\cal J}^{-1}=\Cov{\varrho}
\end{equation}
where the lower bound is attained iff $\varrho=\hat{\varrho}$, which
in turn holds iff $\Pi_2(\Lambda)=0$. This is the case of adaptivity.
The construction of adaptive estimators is a major subject of
semiparametric theory; confer Bickel~(1982; Sections 3 and~4),
Klaassen~(1987), Schick~(1986), and further references mentioned
therein.\qed \end{Rem}
\begin{Rem}\rm [\,bounded influence curves\,]\enskip
The existence of bounded influence curves $\psi\in \Psi$, which may
become relevant for robustness in semiparametric models, proves
equivalent to the following condition
\begin{equation} \label{e.exiShen}
  a' \Lambda \notin \cqlin \partial_2 {\cal Q}
  \qquad \forall\,a\in\R^k,\,a\ne0 \hspace{-1.5em}
\end{equation}
where $\cqlin$ denotes the closed linear span in~$L_1$;
note the difference to~(\ref{e.p.rkJ=k}). The equivalence follows from
Theorem~1 of Shen~(1995) on observing that his condition~\ran$(S')$\ran,
with~$\cqlin ( \partial_2 {\cal Q} + \mbox{constants} )$ in the place
of~$\cqlin \partial_2 {\cal Q}$, because
   $\Ew \Lambda=0$ and $\Ew g=0 \enskip\forall g\in \partial_2 {\cal Q}$,
in fact simplifies to~(\ref{e.exiShen}).
\par     Naturally, 
condition~(\ref{e.exiShen}) is stronger than~(\ref{e.p.rkJ=k}).
When $ \lin \partial_2 {\cal Q} $ has finite dimension, however,
it is closed in both $L_1$ and~$L_2$, and consequentially,
the mere existence of influence curves implies the existence of bounded
ones.\qed \end{Rem}
\begin{Rem}\rm [\,finite dimensions\,] \label{r.p.fidi}
In case $H_{\theta}\subset \R^m$ for some finite dimension~$m$, suppose
the square root densities of model~${\cal Q}$ are $L_2$-differentiable
at~$(\theta_0,\nu_0)$ with respect to the full parameter~$(\theta,\nu)$,
such that~(\ref{e.L2diff}) is satisfied with paths
   $\nu^g_t=\nu_0+t \hspace{4\ei}b$ and $g=b' \Delta$,
where $\Delta\in L_2^m$, $\Ew \Delta=0$, denotes the scores function
for the nuisance parameter~$\nu$ at~$(\theta_0,\nu_0)$.
\par
Then $\partial_2 {\cal Q}= \{\,b' \Delta\mid b\in\R^m\,\}$,
and the Fisher information~${\cal H}$ of model~${\cal Q}$
for the full parameter~$(\theta,\nu)$ at~$(\theta_0,\nu_0)$ is
\begin{equation}
  {\cal H}=  \Cov \pmatrix{ \Lambda \cr \Delta } =
  \pmatrix{ {\cal I} & {\cal C}\cr {\cal C}' & {\cal D}} \qquad
  \mbox{where}\enskip {\cal C}=\Ew \Lambda \,\Delta' \hspace{-1.5em}
\end{equation}
The Fisher information ${\cal D}=\Cov \Delta$ for the parameter~$\nu$
at~$\nu_0$, of model~${\cal Q}_{\theta_0}$, the $\theta_0$-section
of model~${\cal Q}$, is assumed of full rank~$m$.
\par
Then $\Pi_2(\Lambda)= {\cal C}\,{\cal D}^{-1}\Delta$ and
     ${\cal J}= \Cov\bigl(\Lambda- \Pi_2(\Lambda)\big)
              = {\cal I}-{\cal C}\,{\cal D}^{-1}{\cal C}'$.
Moreover, we have ${\cal J}>0$ iff $ {\cal H}>0$,
since $\det {\cal H}=\det {\cal D}\,\det {\cal J}$.
Because ${\cal D}>0$, condition~(\ref{e.exiShen}), too, is 
equivalent to $\rk {\cal H}=k+ \nolinebreak m$. 
\par
In this case,
\begin{equation}\hspace{-1.5em}
   \varrho = {\cal J}^{-1}
   \bigl(\Lambda-{\cal C}\,{\cal D}^{-1}\Delta\bigr) =
   (\EM_k,0_{k \times m})\,
   {\cal H}^{-1} \! \smash{\pmatrix{ \Lambda\cr \Delta}}
   \vphantom{\bigg|}\hspace{-2.5em}
\end{equation}
defines the shortest influence curve---in fact, the first component of
the shortest influence curve for the full parameter, which one usually
is tempted to ascribe to the MLE.
\par
Starting from this function
  $ {\cal H}^{-1}\!\left(\Lambda\atop \Delta\right) $,
bounded influence curves have been constructed explicitly
by HR~(1994), Remark~4.2.11 and 5.5(8),~5.5(9),
if the matrix~$D$ there is specialized to the projection
matrix~$(\EM_k,0_{k \times m})$.\qed   \end{Rem}
Closely related to the orthogonal projection~(\ref{e.cIC}) of influence
curves leading to the canonical influence curve~$\varrho$ is the
Cram\'er--Rao bound for the covariance,
\begin{equation} \label{e.CRao}
  \Cov\psi \ge  {\cal J}^{-1}=\Cov\varrho  \qquad \forall\,\psi\in \Psi
\end{equation}
in the positive definite sense, with equality iff $\psi=\varrho$.
In view of~(\ref{e.ufoN}), this bound concerns the asymptotic covariance
of asymptotically linear estimators. Thus,
the asymptotically linear estimator with canonical influence
curve~$\varrho$ at~$(\theta_0,\nu_0)$ is the asymptotically
most accurate to estimate~$\theta_0$, in model~${\cal Q}$.
\par     That this optimality is not restricted to estimators which are
asymptotically linear, but need to fulfill only a regularity condition
weaker than asymptotic linearity, or may even be arbitrary measurable,
is the subject of the convolution and asymptotic minimax theorems,
respectively; confer, for example,
Bickel et al.~(1993; Theorem~3.3.2), HR~(1994; Theorems 4.3.2,~4.3.4),
van der Vaart~(1998; Theorems 25.20, 25.21, Lemma~25.25).
\begin{Rem}\rm [\,nonlinear projection\,]
These optimality theorems require some structure of the tangent
set~$\partial {\cal Q}$, to be a linear space or at least a convex cone.
In spite of the special structure, the projection in terms of which the
bounds are stated, is generally that on the closed linear
span~$\clin \partial {\cal Q}$.
\par
One exception is the concentration bound for asymptotically median
unbiased estimators by Pfanzagl and Wefelmeyer~(1982; Theorem~9.2.2),
in terms of the projection on a closed convex cone. 
In HR~(2000) we however show that the bound may not possibly be attained,
and  
derive a suitable one-sided bound that is still based on the projection
on~$\cl\partial{\cal Q}$---as opposed to~$\clin \partial {\cal Q}$.\qed
\end{Rem}
\section{The Infinitesimal Robust Setup}
\setcounter{equation}{0}                                      \label{s.B}
In robust statistics, we start with an ideal model
   ${\cal P}= \{\,P_{\theta}\mid \theta\in \Theta\,\}$---from prior
knowledge or nonparametric estimation in advance---which is smoothly
parametrized by some finite ($k$-)dimensional parameter~$\theta$
out of an open subset~$\Theta \subset \R^k$; formally,
  ${\cal P}$~is some model as assumed in Section~\ref{s.P}
but deprived of its nuisance parameter.
Since we do not believe in such a model~${\cal P}$ strictly,
we enlarge its elements~$P_{\theta}$ to certain
neighborhoods~$ U(\theta;r) \subset {\cal M} $ of radius~$r$.
Then the i.i.d.\ observations, under the hypothesis~$\theta$,
may be allowed to follow any law $Q\in \nolinebreak U(\theta;r)$,
while still~$\theta$ has to be estimated.
Thus, the neighborhood model
\begin{equation} \label{e.b.nbdQ}
   {\cal Q}=
   \bigl\{\,Q\bigm| \theta\in \Theta\weg,\,Q\in U(\theta;r)\,\bigr\}
\end{equation}
is obtained, which is clearly semiparametric: For $Q\in U(\theta;r)$,
the deviation $Q-P_{\theta}$ from the ideal~$P_{\theta}$ appears
as nuisance parameter~$\nu$, ranging over the sets of differences
  $H_{\theta}= \{\,Q-P_{\theta}\mid Q\in U(\theta;r)\,\}$,
where $ Q=Q_{\theta, \nu} $ with $\nu\in H_{\theta}$.
In particular, the ideal model~${\cal P}$ is the $\nu_0$-section
of model~${\cal Q}$ at $\nu_0=0$.
\begin{Rem}\rm  [\,nonidentifiability\,]\enskip
This interpretation requires the so-called idealistic robustness approach,
which assumes the existence of an ideal parameter to be estimated even
under deviations from the parametric model.
\par
If one does not start with a true~$\theta$, but seeks~$\theta$ depending
on the real law~$Q$, one runs into the identifiability problem, that is,
multiple solutions~$\theta$ of the equation
  $Q=Q_{\theta, Q-P_{\theta}}=P_{\theta}+Q-P_{\theta}$.
This is the case already for members of the ideal model
  $ Q=P_{\zeta} $ with~$\zeta$ close to~$\theta$
such that $ P_{\zeta}\in \nolinebreak U(\theta;r)$ (if, as usual,
the parametrization is continuous relative to the neighborhoods).
\par
This problem has been dealt with by the `pragmatic' robustness approach,
which defines the parameter by means of functionals that are Fisher
consistent at the ideal model and extend the parametrization to the
neighborhoods. In fact, both approaches lead to the same optimally
robust influence curves and procedures---once the choice of functional is
subjected to robustness criteria;
confer HR~(1994; Preface, Subsection~4.3.3).
So the difference between the two approaches, and hence the difficulty of
the first, seems not essential.
\qed\end{Rem}
\noindent    We specify the neighborhoods~$U(\theta;r)$ to be balls
around~$P_{\theta}$ of radius~$r$ in Hellinger or
total variation distance, or contamination neighborhoods,
\begingroup \mathsurround0em\arraycolsep0em
\begin{eqnarray} \label{e.Uhv}
   U_*(\theta;r) & {}={} & \bigl\{\,Q\in {\cal M}\bigm|
   d_*(Q,P_{\theta})\le r\,\bigr\} \\
\label{e.Uc}
   U_c(\theta;r) & {}={} & \bigl\{\, 
   Q={(1-r)}_{+}\,P_{\theta}+(1\land r)\,M \bigm| 
   M\in {\cal M}\,\bigr\}
\end{eqnarray}\endgroup
where the Hellinger and total variation metrics $d_h$ and~$d_v$
are given by
\begin{equation} \label{e.dhv}
   2 \,d_h^2(Q,P) = \int {\bigl|\sqrt{dQ}-\sqrt{dP}\,\bigr|}^2
   \weg,\qquad    2\,d_v(Q,P) = \int |dQ-dP|
\end{equation}
\noindent    Let us fix $\theta_0\in \Theta$ and $\nu_0=0$, and write~$P$
for the previous $Q=Q_{\theta_0,\nu_0}=\nolinebreak P_{\theta_0}$.
In the sequel, the scores function~$\Lambda$ is that of the ideal
model~${\cal P}$, for~$\theta$ at~$\theta_0$.
\par     Towards the differentiability~(\ref{e.L2diff}) of the neighborhood
model~${\cal Q}_*$ at~$(\theta_0,0)$, depending on the type of
neighborhoods~$U_*(\theta_0;r)$, we introduce the following balls
${\cal G}_* = \nolinebreak {\cal G}_*(\theta_0;r)$
as candidate tangent sets~$\partial_2 {\cal Q}_*$,
\begingroup \mathsurround0em\arraycolsep0em \begin{eqnarray}
\label{e.b.Gh}
   {\cal G}_h & {}={} & \bigl\{\,g\in L_2 \bigm| \Ew g=0\weg,\;
    \Ew g^2\le 8 \hspace{6\ei}r^2\,\bigr\} \\
\label{e.b.Gv}
   {\cal G}_v & {}={} & \bigl\{\,g\in L_2 \bigm| \Ew g=0\weg,\;
    \Ew |g|\le 2 \hspace{6\ei}r\,\bigr\} \\
\label{e.b.Gc}
   {\cal G}_c & {}={} & \bigl\{\,g\in L_2 \bigm| \Ew g=0\weg,\;
    g\ge-r\,\bigr\}
\end{eqnarray}\endgroup
where ${\cal G}_h \subset \sqrt{2}\:{\cal G}_v$ as $ d_v\le \sqrt{2}\,d_h$,
and 
      $ {\cal G}_c \subset {\cal G}_v = {\cal G}_c-{\cal G}_c $
      by~(\ref{e.s.Gc-Gc=Gv}) below.
\par     
The balls~${\cal G}_*$ have already appeared in Bickel~(1981).
\begin{Prop}\sl \label{p.balls}
The tangent sets at~$(\theta_0,0)$ of the neighborhood
model\/~{\rm ${\cal Q}_*$,} for\/~{\rm $*=h,v,c$,} are
\begin{equation}
     \partial_1 {\cal Q}_* = \{\,a' \Lambda \mid a\in\R^k\,\}
\weg, \quad
     \partial_2 {\cal Q}_* = {\cal G}_*
\weg, \quad
     \partial {\cal Q}_* =
     \partial_1 {\cal Q}_* + \partial_2 {\cal Q}_*
\end{equation}
\end{Prop}
\begin{Bew} Invoke bounded approximations~$\Lambda^{(t)}$ of~$\Lambda$
such that $\Ew \Lambda^{(t)}=0$ and, as $t\to0$,
  $ \sup|\Lambda^{(t)}|=\Lo(t^{-1})$ and
  $ \Ew |\Lambda^{(t)}-\Lambda|^2 \to0$.
Given $a\in\R^k$ and any bounded $g\in {\cal G}_*$,
employ the path $\nu_t^g=t \hspace{6\ei}g$ in defining
measures~$Q_t= \nolinebreak Q_{\theta_0+ta,tg}$ by
\begin{equation} \label{e.b.Qt}
   dQ_t= \bigl(1+t \hspace{6\ei}(a' \Lambda^{(t)}+g )\bigr)\,dP
\end{equation}
Then mean square differentiability~(\ref{e.L2diff}) is satisfied, and these
probabilities belong to the neighborhoods~$U_*(\theta_0+ta;tr)$ in the
following, entirely acceptable sense,
\begin{equation} \label{e.b.Qtd*Pt}
   d_*(Q_t,P_{\theta_0+ta}) \le t \hspace{6\ei} r + \Lo(t)
\end{equation}
in the cases $*=h,v$. In the case $*=c$, there exist
approximations~$\widetilde{P}_{\theta_0+ta}$ of~$P_{\theta_0+ta}$,
namely, $\widetilde{P}_{\theta_0+ta}$ with
        $P$~density $1+t_r a' \Lambda^{(t)}$, $t_r=t/(1-tr)$, such that
\begin{equation} \label{e.b.QtinUcPt}
   d_v(\widetilde{P}_{\theta_0+ta}, P_{\theta_0+ta})=\Lo(t)
   \hspace{1.5em}\mbox{and}\hspace{1.5em}
   Q_t\in \widetilde{U}_c(\theta_0+ta;tr)
\samepage\end{equation}
for the contamination balls 
  $\widetilde{U}_c(\theta_0+t \hspace{3\ei} a;t \hspace{3\ei} r)$
  about~$\widetilde{P}_{\theta_0+ta}$.
\par    In either case, we pass to the closure
of~${\cal G}_*\cap L_{\infty}$ in~$L_2$, which is~${\cal G}_*$.
The technical details needed in this proof may be found in HR~(1994):
Remark~4.2.3, Lemma~4.2.4, Lemma~5.3.1, and proof to Theorem~5.4.1\,(a).
\qed\end{Bew}
The tangent sets~${\cal G}_*$ are closed convex, and the
smallest cone and linear space containing either~${\cal G}_*$ is already
the full tangent space $L_2\cap \{\Ew=0\}$, provided only that~$r>0$.
Consequentially, $\Lambda- \Pi_2(\Lambda)=0$ and~${\cal J}=0$
in~(\ref{e.cIC}); in particular, adaptivity fails drastically.
The canonical IC~$\varrho$ is undefined.
\section{The Semiparametric Influence Curve}
\setcounter{equation}{0}                                      \label{s.C}
\noindent    In the robust setup, we therefore modify
definition~(\ref{e.cIC}) of canonical influence curve,
replacing~$\pi_2$ by the nonlinear projection
$\widetilde{\pi}_2 \colon L_2\to \nolinebreak\partial_2 {\cal Q}_* $
on~$\partial_2 {\cal Q}_*= \nolinebreak{\cal G}_*$ itself.
Correspondingly, $\Pi_2$~is replaced by
$ \widetilde{\Pi}_2 ={(\widetilde{\pi}_2,\ldots,\widetilde{\pi}_2)}'
    \colon L_2^k\to (\partial_2 {\cal Q}_*)^k$,
defined coordinatewise.
Thus, the following function~$\tilde{\varrho}_*$,
called semiparametric influence curve, is obtained,
\begingroup \mathsurround0em\arraycolsep0em
\begin{eqnarray} \label{e.cICballs}
   \tilde{\varrho} & {}= {} &
   {\cal K}^{-1}\bigl(\Lambda- \widetilde{\Pi}_2(\Lambda)\bigr)
\hspace{-2em}\\
\noalign{\noindent with scaling matrix\nopagebreak}
\label{e.Kballs}    {\cal K} & {}= {} & \hEw
   \bigl(\Lambda- \widetilde{\Pi}_2(\Lambda)\bigr) \Lambda'
\hspace{-2em}
\end{eqnarray}\endgroup
The definition of~$\tilde{\varrho}$ requires $\det {\cal K}\ne0$.
Rescaling of~$\Lambda- \widetilde{\Pi}_2(\Lambda)$ by~${\cal K}$
ensures Fisher consistency, 
  $\Ew \tilde{\varrho}\hspace{6\ei}\Lambda'= \nolinebreak \EM_k$.
In general
  $ {\cal K}\ne \Cov\bigl(\Lambda- \widetilde{\Pi}_2(\Lambda)\bigr)$,
since residuals are no longer orthogonal to the approximating ball.
\begin{Rem}\rm
The modified projection recipe~(\ref{e.cICballs}),~(\ref{e.Kballs})---%
subtracting from~$\Lambda$ the component explained by the nuisance
parameter, and then rescaling for Fisher consistency---seems no less
plausible than the original one based on linear projection.
Derived only by analogy, the semiparametric influence curve must however
be checked against a mathematical solution to some suitable extension
of the Cram\'er--Rao bound, or convolution and asymptotic minimax theorems,
in the semiparametric/robust 
setup with full tangent balls.\qed\end{Rem}
\noindent    The following approximation lemma is well-known
and will be applied to the balls $G=\nolinebreak{\cal G}_*$,
the space $X=\nolinebreak L_2$, and the coordinates~$x$
of~$\Lambda$; then $\tilde{g}=\widetilde{\pi}_2(\Lambda_j)$.
\begin{Lem}\sl \label{l.c.approx}
Let~$G$ be a nonempty closed and convex subset of some Hilbert
space~{\rm $X$,} and~$x\in X$.
Then the minimum norm problem
\begin{equation}
   {|x-g|}^2=\min {!}\qquad g\in G
\end{equation}
has a unique solution~{\rm $\tilde{g}\in G$,} which is characterized by
\begin{equation} \label{e.gmin}
  \langle x- \tilde{g}| g - \tilde{g} \rangle \le 0
  \qquad \forall\,g\in G
\end{equation}
\end{Lem}
\noindent    
In the sequel, ${\cal I}=\Cov \Lambda=({\cal I}_{i,j})$
and $ \hat{\varrho} = {\cal I}^{-1}\Lambda $ denote Fisher information
(of full rank~$k$) and the canonical influence curve,
of the ideal model~${\cal P}$ at~$\theta_0$.
\par     We now determine the semiparametric influence curves
$\tilde{\varrho}_h$, $\tilde{\varrho}_v$, $\tilde{\varrho}_c$ for the
Hellinger, total variation, and contamination neighborhood models,
respectively.
\begin{Thm}{\rm [\,Hellinger model\,]}\enskip      \sl \label{t.rhoh}
The semiparametric IC~$\tilde{\varrho}_h$ exists iff
\begingroup \mathsurround0em\arraycolsep0em
\begin{eqnarray} \label{e.rh}
  8 \,r^2 & {} < {} &
  \min\nolimits_{j=1,\ldots,k}\, {\cal I}_{j,j} \\ 
\noalign{\noindent And then} \label{e.rhoh} 
  \tilde{\varrho}_h & {} = {} & \hat{\varrho} = {\cal I}^{-1}\Lambda
\end{eqnarray}\endgroup
\end{Thm}
\begin{Bew} In the case $k=1$ we have $\widetilde{\pi}_2=\gamma \Lambda$
with $\gamma=$~positive root of the minimum of~$1$
and~$8 \,r^2 \!/{\cal I}$.
Indeed, by Cauchy--Schwarz, for every $g\in {\cal G}_h$,
\begin{equation}
  \langle \Lambda- \gamma \Lambda|g\rangle =
  (1 - \gamma)\langle \Lambda|g \rangle \le
  (1 - \gamma) \sqrt{8}\,r\,{\cal I}^{1/2} =
  (1 - \gamma) \hskip6\ei\gamma \,{\cal I} =
  \langle \Lambda- \gamma \Lambda| \gamma \Lambda\rangle
\end{equation}
For general $k\ge1$, this implies that
  $\Lambda- \widetilde{\Pi}_2(\Lambda)= D \Lambda$ and
  ${\cal K}= D \,{\cal I}$ with matrix $D=\diag(1- \gamma_j)$, where
  $0\le\gamma_j\le1$, and $\gamma_j=1 $ iff
  $ {\cal I}_{j,j}\le 8 \hskip6\ei r^2$.
\qed\end{Bew}
\begin{Thm}{\rm [\,total variation\,]}\enskip         \sl \label{t.rhov}
The semiparametric IC~$\tilde{\varrho}_v$ exists only if
\begingroup \mathsurround0em\arraycolsep0em
\begin{eqnarray}
\label{e.rv} 2 \,r &{}<{}& \min \nolimits_{j=1,\ldots,k}\,\Ew |\Lambda_j|\\
\noalign{\vspace{.66\belowdisplayshortskip}\pagebreak[0]
\mathsurround=\msu\noindent
  And then $\tilde{\Lambda}^{(v)}=\Lambda- \widetilde{\Pi}_2(\Lambda)$
  has coordinates
\nopagebreak\vspace{.66\abovedisplayskip}}
\label{e.rhov}
  \tilde{\Lambda}_j^{(v)} & {}={} & v'_j \lor \Lambda_j \land v''_j \\
\noalign{\vspace{\belowdisplayshortskip}\pagebreak[0]
\mathsurround=\msu\noindent
  where the clipping constants $ v'_j  < 0 < v''_j $
  are uniquely determined by \nopagebreak\vspace{\abovedisplayskip}}
\label{e.gv}   \hEw ( v'_j - \Lambda_j)_{+} & {}={} & r =
   \hEw ( \Lambda_j - v''_j \hspace{6\ei})_{+}
\end{eqnarray}\endgroup
\end{Thm}
\begin{Bew} Obviously, $\Lambda_j-\widetilde{\pi}_2(\Lambda_j)=0$
iff $\Ew|\Lambda_j|\le 2 \hspace{6\ei}r$. Thus assume~(\ref{e.rv}).
\par     In case $k=1$, in order to minimize $\Ew(\Lambda-g)^2$
for $g\in {\cal G}_v$, we set up a Lagrangian
$  
   \hEw \bigl( (\Lambda-g)^2
   + 2 \hspace{6\ei} \alpha \hspace{6\ei} g
   + 2 \hspace{6\ei} \beta \hspace{6\ei} |g| \,\bigr)
$  with some unspecified real multipliers,
and try to minimize the integrand
$ I(g)= (\Lambda-g)^2
   + 2 \hspace{6\ei} \alpha \hspace{6\ei} g
   + 2 \hspace{6\ei} \beta \hspace{6\ei} |g|
$  at each point.
\par     A minimizing value~$\tilde{g}=0$ means that
$ \Lambda^2 \le (\Lambda-g)^2
   + 2 \hspace{6\ei} \alpha \hspace{6\ei} g
   + 2 \hspace{6\ei} \beta \hspace{6\ei} g
$  for all numbers $g>0$; that is, $\Lambda- \alpha\le \beta$,
and
$ \Lambda^2 \le (\Lambda-g)^2
   + 2 \hspace{6\ei} \alpha \hspace{6\ei} g
   - 2 \hspace{6\ei} \beta \hspace{6\ei} g
$  for all numbers $g<0$; that is, $\Lambda- \alpha\ge -\beta$.
This is the case when $\Lambda- \tilde{g}= \Lambda$.
\par     If $\tilde{g}>0$, then the derivative
$d I (\tilde{g})=0$ gives 
$\Lambda- \tilde{g}= \alpha+ \beta$.
If $\tilde{g}<0$, $dI (\tilde{g})= \nolinebreak 0$
gives $\Lambda- \tilde{g}= \alpha-\beta$.
These are the cases when $\Lambda - \alpha > \beta$,
respectively when $\Lambda - \alpha < -\beta$.
\par     Altogether,
$ \Lambda- \tilde{g} =
  {(-\beta)\lor(\Lambda- \alpha)\land \beta} \: + \alpha
  = (\alpha- \beta)\lor \Lambda \land(\alpha+ \beta) $ seems to be
the necessary form of $\tilde{q}= \Lambda- \tilde{g}$.
\par     Now define $\tilde{q}= v'\lor \Lambda\land v''$
by means of the unique solutions $v'<0<v''$ of
$ \Ew(v'- \Lambda)_+=r=\Ew(\Lambda-v'' \hspace{3\ei})_+ $,
which is a matter of continuity (dominated convergence theorem),
monotony (strict), and the intermediate value theorem.
We shall verify that this~$\tilde{q}$ minimizes~$\Ew q^2$
subject to $\Ew q=0$, $\Ew|\Lambda-q|\le 2 \hspace{6\ei}r$.
\par     By the definition of~$\tilde{q}$,
$  \Ew (\Lambda-q)\tilde{q}
   \le v''\Ew (\Lambda-q)_+ - v'\Ew (q - \Lambda)_+
$,  which is less or equal
$   r \hspace{6\ei}(v'' - v' \hspace{3\ei})
    = \Ew (\Lambda- \tilde{q})\tilde{q} $.
Thus $\Ew (- \tilde{q})(q- \tilde{q})\le0$, which is~(\ref{e.gmin}).%
\qed\end{Bew}
\begin{Thm}{\rm[\,contamination\,]}\enskip           \sl \label{t.rhoc}
The semiparametric IC~$\tilde{\varrho}_c$ exists only if
\begingroup \mathsurround0em\arraycolsep0em
\begin{eqnarray} \label{e.rc}
  r & {}< {} &  - \max \nolimits_{j=1,\ldots,k} \,\winf_{P} \Lambda_j
\hspace{-2em}\\ \noalign{\vspace{\belowdisplayshortskip}\pagebreak[0]
\mathsurround=\msu\noindent   where
   $\winf_P$ denotes the $P$~essential infimum. And then
   $\tilde{\Lambda}^{(c)}=\Lambda- \widetilde{\Pi}_2(\Lambda)$
   has coordinates\nopagebreak\vspace{\abovedisplayshortskip}}
\label{e.rhoc}
   \tilde{\Lambda}_j^{(c)} & {}={} & (\Lambda_j + r) \land u_j \\
\noalign{\vspace{\belowdisplayshortskip}\pagebreak[0]
\mathsurround=\msu\noindent   with clipping constant $ u_j>0$
   uniquely determined by\nopagebreak\vspace{\abovedisplayskip}}
\label{e.gc} 0 & {}= {} & \hEw {(\Lambda_j + r)\land u_j}
\end{eqnarray}\endgroup \end{Thm}
\begin{Bew} Obviously, $\Lambda_j-\widetilde{\pi}_2(\Lambda_j)=0$
iff $\Lambda_j\ge-r$~a.e.$P$. Thus assume~(\ref{e.rc}).
\par     In case $k=1$, in order to minimize $\Ew(\Lambda-g)^2$
for $g\in {\cal G}_c$, we pass to the equivalent problem of
minimizing~$\Ew q^2$ subject to $\Ew q=0$, $q\le \Lambda+r$,
for which we minimize a Lagrangian
$\Ew q^2- 2 \hspace{6\ei}u\Ew q =
 \Ew (q - u)^2+ \mbox{constant}$, subject to $q\le \Lambda+r$.
Doing this pointwise, the necessary form seems
$\tilde{q}=(\Lambda+r)\land u$.
\par     Now consider the function $f(s)=\Ew(\Lambda+r)\land s$
for $s\ge0$. It is monotone, continuous [\,dominated convergence
applies since $-(\Lambda+r)_{-}\le f\le (\Lambda+r)_{+}$], and
has limits $-\Ew(\Lambda+r)_{-}<0$ and~$r\ge0$ at~$0$ and~$\infty$,
respectively. Thus~$f$ has a zero~$u>0$, which we use to define
$\tilde{q}=(\Lambda+r)\land u$.
(Only in case~$r=0$, may~$u$ be nonunique, but then
$\tilde{q}=\Lambda$.) By construction, $\tilde{q}$~satisfies
the side conditions $\Ew q=0$, $q\le \Lambda+r$.
\par     To prove~$\tilde{q}$ optimal, let~$q\in L_2$ be any such function.
Then $q\le \tilde{q}= \nolinebreak\Lambda+r$ as soon as $\tilde{q}< u$.
Thus $(u - \tilde{q})(u- q)$ is always greater
or equal to $(u- \tilde{q})^2$. Consequentially,
$\Ew(- \tilde{q})(q- \tilde{q})= 
 \Ew(u- \tilde{q})(q- u+ u-\tilde{q})\le0$; 
which is~(\ref{e.gmin}).%
\qed\end{Bew}
\begin{Rem}\rm  \label{r.K>0}
In Theorems~\ref{t.rhov} and~\ref{t.rhoc}, conditions~(\ref{e.rv})
and~(\ref{e.rc}), respectively, ensure that
   $ \hEw \Lambda_j^{(*)}\Lambda_j >0 $ for $*=v,c$.
This may be seen by writing
\begingroup \mathsurround0em\arraycolsep0em \begin{eqnarray}
   \hEw \Lambda_j^{(v)}\Lambda_j &{}={}& \hEw \Lambda_j^{(v)}
   \Lambda_j^{(v)} + r \hspace{6\ei} (v_j'' - v_j' \hspace{3\ei}) \\
\noalign{\vspace{\belowdisplayshortskip}\pagebreak[0]
\mathsurround=\msu\noindent
   where $ r \hspace{6\ei} (v_j'' - v_j' \hspace{3\ei}) > 0 $
   unless $r=0$ (and then $\Lambda_j^{(v)}= \Lambda_j$,
   and~${\cal I}_{j,j}>0$),  respectively by writing
\mathsurround=\msu\noindent   \nopagebreak\vspace{\abovedisplayskip}}
   \hEw \Lambda_j^{(c)}\Lambda_j &{}={}& \hEw \Lambda_j^{(c)}\Lambda_j^{(c)}
   + \hEw \Lambda_j^{(c)} \bigl( \Lambda_j + r - \Lambda_j^{(c)}\, \bigr)
\end{eqnarray}\endgroup
where $\Lambda_j+r\le u_j$~a.e.$P$ only if $r=0$ (and again
   $\Lambda_j^{(c)}= \Lambda_j$,~${\cal I}_{j,j}>0$).
\par     However, whether condition~(\ref{e.rv}),
respectively~(\ref{e.rc}), for dimension $k>\nolinebreak 1$ already
imply the nonsingularity of~${\cal K}$, 
hence the existence of~$\tilde{\varrho}_v$, respectively
of~$\tilde{\varrho}_c$, is unclear.\qed
\end{Rem}
\begin{Rem}\rm  The optimization problems 
of this section resemble those that 
determine robust influence curves, 
however with three distinctions:
\begin{ABCa}\item
         The approximation $\Ew|\Lambda-g|^2=\min{!}$,
         instead of $\Ew |\psi|^2=\min{!}$.
\item    The $L_2$,$L_1$, and $\winf_P$~bounds on tangents 
         translate into bounds on influence curves 
         in the dual norm~$\sup_{g\in {\cal G}_*}|{\Ew{\psi g}}|$,
         for $*=h,v,c$, respectively.
\item    There is no condition on tangents that would correspond to the
         Fisher consistency $\Ew \psi \Lambda' = \nolinebreak \EM_k $
         of influence curves. 
\qed \end{ABCa}
\end{Rem}
Thus, at least, certain features of robust influence curves are recovered.
\section{Comparison of Semiparametric and Robust Estimators}
\setcounter{equation}{0}\label{s.R}
More precisly, the semiparametric
  recipe~(\ref{e.cICballs}),~(\ref{e.Kballs})
will be judged under a certain etimator risk.
How does the semiparametric estimator---the asymptotically linear estimator
with semiparametric influence curve~$\tilde{\varrho}_*$---compare with the
robust estimator---the asymptotically linear estimator with robust influence
curve~$\eta_*$ that, by definition, minimizes maximum asymptotic mean square
error of asymptotically linear estimators? The maximum is evaluated over
shrinking neighborhoods~$U_*(\theta_0; r /\!\sqrt{n}\,)$, as the sample
size~$n$ tends to infinity, with starting radius~$r\ge0$---henceforth,
radius~$r$---fixed. For asymptotically linear estimators, this maximum
asymptotic MSE naturally extends the covariance criterion employed in
the Cram\'er--Rao bound to the infinitesimal robust setup.
\begin{Rem}\rm
An extension of asymptotic maximum~MSE over neighborhoods,
from asymptotically linear to arbitrary estimators $S=(S_n)$,
employing a risk such as
\begin{equation}  \label{e.rAM}
  \lim_{b\to \infty} \lim_{c\to \infty} \limsup_{n\to \infty}
  \hspace{12\ei}\sup_{|t|\le c}
  \hspace{18\ei} \sup_{Q\in U_n(t; r \hskip3\ei)}
  \hspace{18\ei}   \int b\land |R_n|^2   \,dQ^n
\end{equation}
where $ U_n(t; r \hskip6\ei) =
      U_* (\theta_0+t/\!\sqrt{n}\,, r /\!\sqrt{n}\,) $
of fixed radius~$ r $, and $R_n=\sqrt{n}\,(S_n- \theta_0)$,
has not been achieved.
Theorem~4.1(A) of HR~(1981\,b), which admits arbitrary estimators,
is restricted to one sided confidence probabilities,
dimension~$k= \nolinebreak 1$, and total variation, contamination
neighborhoods (for which least favorable probability pairs exist).
Therefore, except in this special case, 
the comparison of semiparametric and robust ICs is bound to
asymptotically linear estimators.\qed\end{Rem}  
For the estimation of~$\theta_0$,
over shrinking neighborhoods $U_*(\theta_0;  r /\!\sqrt{n}\,)$,
radius~$ r $, we consider a weighted~MSE with nonnegative
bias weight~$\beta$.
In the case of estimators of~$\theta_0$ that are asymptotically linear
with influence curves~$\psi$ at~${\theta_0}$,
the maximum asymptotic weighted mean square error is
\begin{equation} \label{e.mMSE}
   \MSE_{*}(\psi;\beta,r)=
   \Ew|\psi|^2 +  \beta \hspace{6\ei}r^2 \omega_{*}^2(\psi)
\end{equation}
As for the derivation of this risk with weight $\beta=1$,
the bias terms~$\omega_{*}(\psi)$, and the minimization
of~$\MSE_{*}(\psi;\beta,r)$ for $\psi\in\Psi$, which determines
the robust influence curve~$\eta_*$ uniquely, 
confer HR~(1994; Subsection~5.5.2).
\par     The influence curves $\Psi= \Psi_{\theta_0}$,
and asymptotic linearity of estimators, are defined with respect to
the ideal model~${\cal P}$ at~$\theta_0$.
\subsection{Coincidence in Hellinger Model}
Hellinger bias, according to HR~(1994; Proposition~5.5.3),
is given in terms of the maximum eigenvalue of the covariance,
$ \omega^2_h(\psi) = 8 \hspace{6\ei}\maxev \Cov\psi $.
In view of the Cram\'er--Rao bound~(\ref{e.CRao}), therefore,
Hellinger risk~$\MSE_{h}(\ldotp;\beta,r)$ is minimized by the canonical
IC~(\ref{e.rhoP}): $ \hat{\varrho} =  {\cal I}^{-1}\Lambda $, for every
  $ \beta, r\in [\,0,\infty) $.
Theorem~\ref{t.rhoh} thus yields the following coincidence. 
\begin{Thm}\sl \label{p.etHro}
Assume\/~{\rm (\ref{e.rh}):
   $ 8 \,r^2 <  \min_{j=1,\ldots,k}\, {\cal I}_{j,j} $.} 
Then the semiparametric IC~$\tilde{\varrho}_h$ agrees with the
robust IC~$\eta_h$,
\begin{equation} \label{e.rhoHeta}
   \tilde{\varrho}_h = \hat{\varrho} = {\cal I}^{-1}\Lambda = \eta_h
\end{equation}
minimizing\/~{\rm $\MSE_h(\ldotp;\beta,r)$,}  
for every $\beta\in[\,0,\infty)$.
\end{Thm}
In principle, the coincidence is a first justification of
the semiparametric recipe.
The value of this result, however, is somewhat diminished since
Hellinger balls, in certain respects, are deemed too small;
confer Bickel~(1981; Th\'eor\`eme~8) and HR~(1994; Example~6.1.1).
The gross error neighborhoods (total variation, contamination) seem
in practice more suitable for robustness.
\begin{Rem}\rm  \label{r.r.Hopt} Identity~(\ref{e.rhoHeta}) implies
equality $\Cov\tilde{\varrho}_h=\Cov\hat{\varrho}$ in~(\ref{e.Cadapt}),
with the semiparametric and robust IC~$\hat{\varrho}=\eta_h$ in the place
of the canonical IC, which might suggest adaptivity. 
However, due to bias, covariance alone does not define the right
risk in the Hellinger model~${\cal Q}_h$, which is why~$\MSE_h$
is used. Clearly,  
\begin{equation} \label{e.MSEh0}
  \MSE_h(\eta_h;\beta,r) =
  \tr {\cal I}^{-1}+8 \hskip6\ei \beta \hspace{6\ei}
  r^2 \maxev {\cal I}^{-1}
  > \tr {\cal I}^{-1} = \MSE_h(\hat{\varrho};\beta,0)
\end{equation}
if only $\beta \hspace{6\ei}r>0$.
Thus, despite~$\hat{\varrho}$ achieves minimum~$\MSE$ in model~${\cal Q}_h$
as well as in~${\cal P}$, strict inequality holds in~(\ref{e.MSEh0}),
so adaptivity is violated; Hellinger neighborhoods do not go for free.
\qed\end{Rem}
\subsection{Relations for Total Variation}
\subsubsection{\protect\boldmath Dimension $k=\bf1$}
\begin{figure}                      
\hfill\includegraphics*[trim= 13 230 130 90, %
                              height=2in,width=3in]{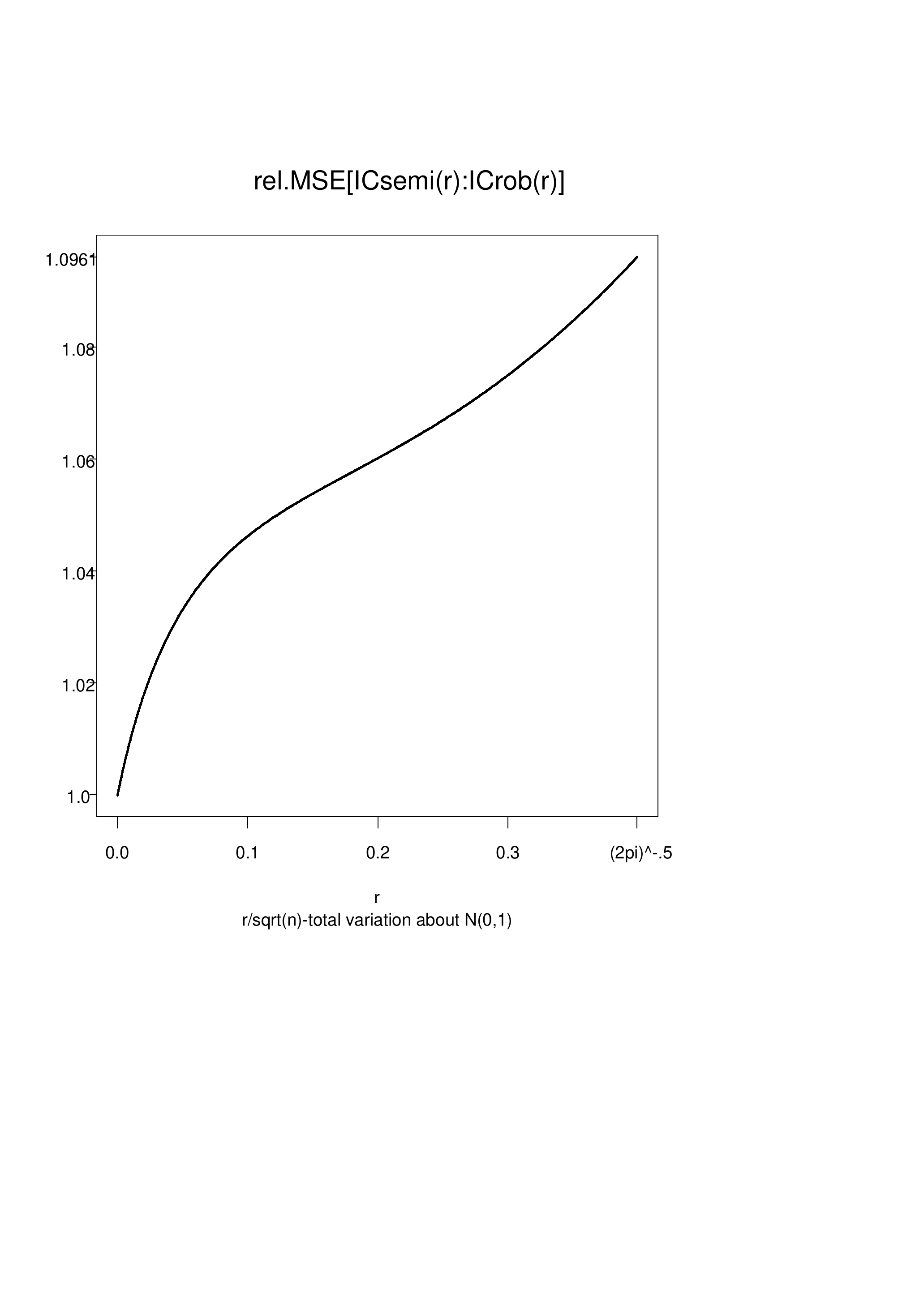}\hfill\mbox{} 
\caption{\protect\small   $\relMSE$ of semiparametric
   IC~$\tilde{\varrho}_v$ vs.~robust IC~$\eta_v$,
   for $0\le r < 1/ \!\protect\sqrt{2 \pi}\,$.}
\label{f.relMSEr}\end{figure}

Total variation bias in one dimension,
according to HR~(1994; Proposition~5.5.3),
is $ \omega_v(\psi)=\wsup_P\psi-\winf_P\psi$.
The robust IC~$\eta_v$ minimizing~$\MSE_v(\ldotp;\beta,r)$
is given by HR~(1994; Theorem~5.5.7), with~$\beta \hspace{6\ei} r^2$
replacing~$\beta$ there. Thus,%
\nopagebreak\begingroup \mathsurround0em\arraycolsep0em
\begin{eqnarray}
\label{e.etav1}  \eta_v &{}={}& c'\lor A\Lambda\land c'' \\
\noalign{\vspace{\belowdisplayshortskip}\pagebreak[0]
\mathsurround=\msu\noindent
  for any numers $c'<0<c''$ and~$A$ such that
  $\Ew \eta_v=0$, $\Ew \eta_v \Lambda=1$, and
\nopagebreak\vspace{\abovedisplayskip}}
\label{e.etav1b}  \hspace{-1.66em}
  \beta \hspace{6\ei}r^2 (c''-c') &{}={}& \hEw (c'-A \Lambda)_{+}
\end{eqnarray}\endgroup
The following result justifies the semiparametric
recipe~(\ref{e.cICballs})--(\ref{e.Kballs}) if one accepts
the particular bias weight implicitly defined by~(\ref{e.betav}).
\begin{Thm}\sl \label{t.retv1}  
Assume\/~{\rm (\ref{e.rv}): $ r < \Ew \Lambda_+ $.}
Then the semiparametric IC~$\tilde{\varrho}_v$ agrees with the robust
IC~$\eta_v$ minimizing {\rm $\MSE_v(\ldotp;\beta,r)$,} 
iff bias weight~$\beta=\beta(r)$ is chosen such that
\begin{equation}\label{e.betav}
   \beta^{-1} = r \hskip6\ei (v''- v' \hskip6\ei)
\hspace{-1em}\end{equation}
where $v'=v'(r)<0<v''(r)=v''$ are determined by
$$ \hEw ( v' - \Lambda)_{+} = r =
   \hEw ( \Lambda - v'' \hspace{6\ei})_{+}  \hspace{-18\ei}
   \eqno (\ref{e.gv})\hspace{\parindent} $$
\end{Thm}
\begin{Bew} Theorem~\ref{t.rhov} supplies
  $\tilde{\varrho}_v= A\, v'\lor \Lambda\land v''$ with
clipping constants $v', v''$ determined by~(\ref{e.gv})
and rescaling constant $A^{-1}={\cal K}>0$ (Remark~\ref{r.K>0}).
\par     Thus $\tilde{\varrho}_v$ attains form~(\ref{e.etav1})
with $c'=v'A$ and $c''=v''A$; in particular,
$ \beta \hspace{6\ei}r^2 (c''-c') = \nolinebreak
  \beta \hspace{6\ei}r^2 (v''- v' \hskip6\ei) A $.
Since $ Ar = A \hEw (v'-\Lambda)_{+} = \hEw (c'-A \Lambda)_{+} $
by~(\ref{e.gv}), equation~(\ref{e.etav1b}) is the same as~(\ref{e.betav}).
\qed\end{Bew}
Bias weight $\beta=1$, in view of~(\ref{e.rAM}), seems the most
natural choice. Then the semiparametric IC~$\tilde{\varrho}_v$
minimizes~$\MSE_v(\ldotp;1,r_1)$, since it is the robust
IC~$\eta_v$ for this radius~$r_1$, iff
\begin{equation} \label{e.rbv1}
    r_1^{-1}=v''(r_1)- v' (r_1)
\end{equation}
Let us keep bias weight $\beta=1$.
Then the semiparametric IC~$\tilde{\varrho}_v$ defined for radius~$r$
minimizes the risk $\MSE_v(\ldotp;1,R(r))$ for another radius~$R(r)$
given by
\begin{equation} \label{e.Rr}
   R^2(r)= r \big/\bigl(v''(r)- v' (r)\bigr) =  r^2 \beta(r)
\end{equation}
since~$\tilde{\varrho}_v$ is of form (\ref{e.etav1}) and~(\ref{e.etav1b}),
hence is the robust~$\eta_v$, for this radius~$R(r)$.
\begin{figure}
\vspace{\bigskipamount} 
\hspace{-1em}\includegraphics[trim=20 220 20 220,
                       height=2.5in,width=2.5in]{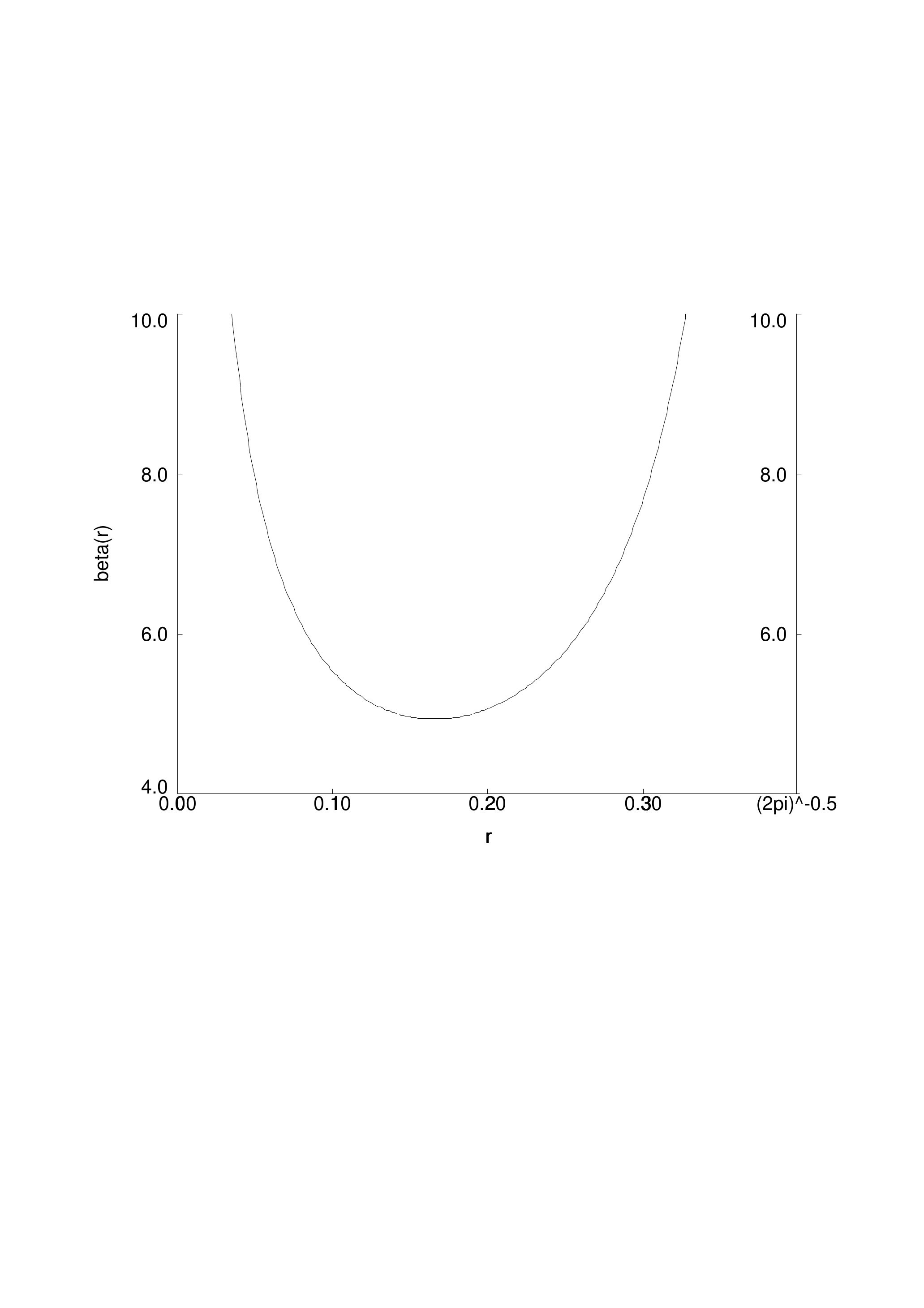} 
\includegraphics[trim=20 220 20 220,
                       height=2.5in,width=2.5in]{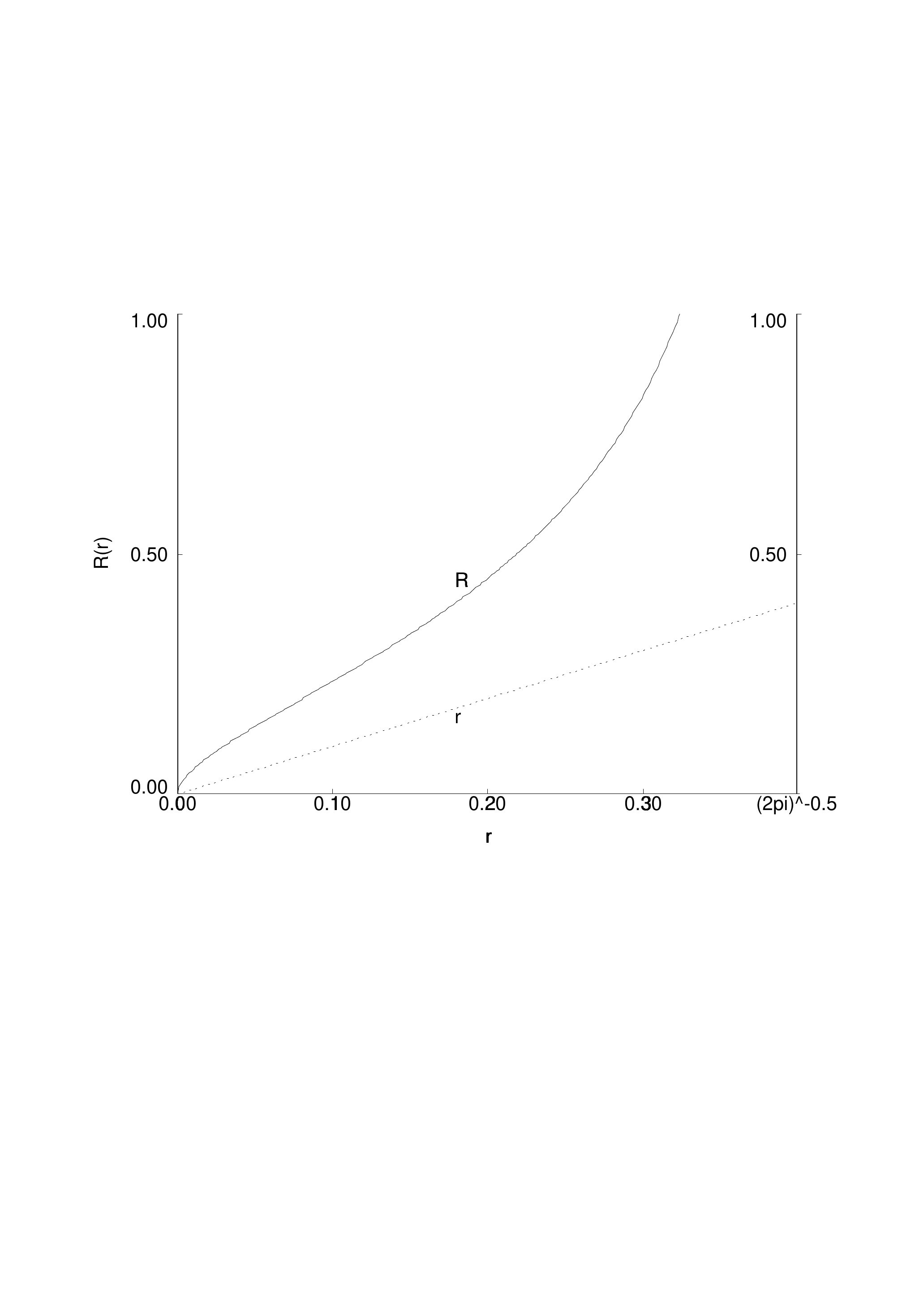}
 \vspace{-3\bigskipamount}\caption{\protect\small
         Bias weight~$ \beta(r)$ 
         and radius~$R(r)$ 
         versus radius $0< r < 1/ \!\protect\sqrt{2 \pi}\,$,
for total variation neighborhoods~$U_v(\theta,r/\!\protect\sqrt{n}\,)$
about the ideal location model $P_{\theta}= {\cal N}(\theta,1)$.}
\label{f.bRr}\end{figure}
\begin{Exa}\rm  For the standard normal location model
  $P_{\theta}= {\cal N}(\theta,1)$, Figure~\ref{f.bRr} shows the bias
weight~$\beta(r)$ and the radius~$R(r)$ defined by~(\ref{e.betav})
and~(\ref{e.Rr}), respectively.
The function~$\beta(\ldotp)$ has singularities at~$0$ and the right
boundary, which is~$1/\!\sqrt{2 \pi}\,=0.3989$, and attains its
minimum value $\beta_{\rm min}=4.8662$ at~$r_{\rm min}=0.1668$.
In particular, no radius~$r_1$ for which $\beta(r_1)=1$ exists.
\par
The radius~$R(r)$ descends to~$0$, hence $\beta(r)=\Lo(r^{-2})$, 
as $r\to0$, and rises to~$\infty$ as $r\to1/\!\sqrt{2 \pi}\,$.
Since~$ R(r)/r = \sqrt{\beta(r)}\,$ is always larger
than~$\sqrt{\beta_{{\rm min}}}\,$, the semiparametric
IC~$\tilde{\varrho}_v$ safeguards against more than double
the amount of contamination assumed in its
definition~(\ref{e.cICballs})--(\ref{e.Kballs}) and,
as~$\beta(r)>\beta_{\rm min}$, is typically even more pessimistic.
\par     Nevertheless, the efficiency loss incurred by the
semiparametric~IC is not dramatic.
Figure~\ref{f.relMSEr} plots the relative maximum asymptotic MSE,
\begin{equation} \label{e.relMSEvr}
         \relMSE(\eta_v,r):=
         \frac{\MSE_v(\eta_v;1,r)}{\MSE_v(\tilde{\varrho}_v;1,r)}
  \weg, \qquad 0\le r < \frac{1}{\sqrt{2 \pi}\,} \hspace{-2em}
\end{equation}
of the semiparametric IC $\tilde{\varrho}_v=\tilde{\varrho}_{v,r}$
in comparison to the robust IC $\eta_v=\eta_{v,r}$, as a function
of the radius~$r$. 
   $\relMSE(\eta_v,r)$ smoothly increases from its
minimum value~$1$ at $r=0$ to its supremum~$1.0961$
as $r\to1/\!\sqrt{2 \pi}\,$.
\par
Therefore, the semiparametric~IC never needs more than $10\%$
additional observations, in order to achieve the same accuracy
in terms of~$\MSE_v(\ldotp;1,r)$ as the robust~IC.
\qed \end{Exa}
\subsubsection{Confidence risk}
The asymptotic maximum risk considered in HR~(1981\,b), instead of
mean square error, and bounded from below for arbitrary estimators~$(S_n)$,
is based on right and left confidence probabilities as follows,
\begin{equation} \label{e.rcAM}
  \lim_{c\to \infty} \limsup_{n\to \infty} \hspace{12\ei}
  \sup_{|t|\le c} \hspace{18\ei} \sup_{Q\in U_n(t; r \hskip3\ei)}
  {Q^n (R_n<-\tau)} \lor {Q^n (R_n>\tau)}
\end{equation}
where $ U_n(t; r \hskip6\ei) =
        U_v (\theta_0+t/\!\sqrt{n}\,, r /\!\sqrt{n}\,) $
of fixed radius~$ r $, and $\tau\in (0,\infty)$ is some interval half-width.
As already in~(\ref{e.rAM}), the standardization
$R_n=\sqrt{n}\,(S_n- \nolinebreak\theta_0)$ is needed only
for the description of the asymptotic minimax estimator
as an asymptotically linear one.
\begin{Thm}\sl  \label{t.CIrhov}
Assume~{\rm (\ref{e.rv}):} $r<\Ew \Lambda_+$. 
Then the semiparametric IC~$\tilde{\varrho}_v$ agrees with the
robust IC~$\eta_v$ with respect to confidence risk~{\rm (\ref{e.rcAM})}
iff we choose half-width 
\begin{equation} \label{e.tau1}
   \tau=\tau(r)=1
\end{equation}
\end{Thm}
\begin{Bew}
According to HR~(1981\,b; Theorems~4.1(A)--4.3; 1980; Theorem~3.1),
for radius
\begin{equation} \label{e.tau/rcAM}
  r< \tau \Ew \Lambda_+
\end{equation}
the estimator~$(S_n)$ minimizing risk~(\ref{e.rcAM}) 
is asymptotically linear at~$\theta_0$ with IC~$\eta_v$ of form
(\ref{e.rhov}) and~(\ref{e.gv}), however, with~$r$ in~(\ref{e.gv})
replaced by~$r/\tau$.
\par    Thus, the semiparametric IC~$\tilde{\varrho}_v$ is the robust
IC~$\eta_v$ iff $\tau=1$ in risk~(\ref{e.rcAM}).
And then, condition~(\ref{e.tau/rcAM}) on~$r$ 
is the same 
as~(\ref{e.rv}).\qed \end{Bew}
\subsubsection{\protect\boldmath Dimension $k>\bf1$}
Exact total variation bias for more than one dimension is rather unwieldy,
$ \omega_v(\psi) = \sup_{|e|=1}\wsup_P{e'\psi} - \winf_P{e'\psi} $,
where $\sup_{|e|=1}$ extends over all unit vectors in~$\R^k$; confer
HR~(1994; Proposition~5.3.3).
Approximate versions $\omega_{v;2}^2(\psi)$ and $\omega_{v;\infty}(\psi)$
have been defined by the Euclidean and $\sup$~norms in~$\R^k$ of the vector
of coordinate biasses~$ \omega_v(\psi_j) $, respectively, which bound the
exact bias from below and above:
$ \omega_{v;\infty} \le \omega_{v} \le \omega_{v;2} \le
  \sqrt{k}\:\omega_{v;\infty} $.
According to HR~(1994; Theorems~5.5.6--7) on one hand, the robust
ICs~$\eta_v$ minimizing either risk $\MSE_{v;s}(\ldotp;\beta,r)$
have the coordinates
\begin{equation} \label{e.etvj}
   \eta_j= c'_j\lor A_j \Lambda \land c''_j
\end{equation}
with any numbers $c'_j<0<c''_j$ and row vectors $A_j\in\R^k$
such that the side conditions $\Ew \eta_v=0$ and
$\Ew \eta_v \Lambda'=\EM_k$ are met.
Moreover, the clipping constants satisfy
\begin{equation}
   \beta \hspace{6\ei}r^2 (c''_j-c'_j)= \hEw {(c'_j-A_j \Lambda)}_{+}
\end{equation}
in case $s=2$, whereas, in case $s=\infty$,
   the differences $c''_j-c'_j$ are all the same
\begin{equation}
   \beta \hspace{6\ei}r^2 (c''_j-c'_j) = \hEw {(c'_1-A_1 \Lambda)}_{+}
   + \cdots + \hEw {(c'_k-A_k \Lambda)}_{+}
\end{equation}
By Theorem~\ref{t.rhov} on the other hand, with clipping constants
  $v'_j<0<v''_j$ defined by~(\ref{e.gv}), and $(A^{j,i})^{-1}= {\cal K}$
by~(\ref{e.Kballs}), the semiparametric IC~$\tilde{\varrho}_v$
has coordinates
\begin{equation} \label{e.r.spICv}
   \tilde{\varrho}_j= A^{j,1}\, v'_1 \lor \Lambda_1 \land v''_1
   + \cdots + A^{j,k}\, v'_k \lor \Lambda_k \land v''_k
\end{equation}
Thus, the order of clipping and linear combination is interchanged
in $\tilde{\varrho}_v$ and~$\eta_v$. So $\tilde{\varrho}_v$ resembles,
but does not exactly match, the robust~$\eta_v$, therefore does not
minimize either risk $\MSE_{v;s}(\ldotp;\beta,r)$, $s=2,\infty$,
if only $\beta \hspace{6\ei}r>0$.
\par   However, the bias terms $\omega_{v;s}$ are only bounds for the
exact bias~$\omega_{v}$, while~$\tilde{\varrho}_v$ ought to be compared
with the minimizer of the exact risk~$\MSE_v(\ldotp;\beta,r)$.
And, at least, $\tilde{\varrho}_v$ has finite biasses
$\omega_{v;s}(\tilde{\varrho}_v)$ and~$\omega_{v}(\tilde{\varrho}_v)$,
hence finite risks $\MSE_{v;s}(\tilde{\varrho}_v;\beta,r)$,
and~$\MSE_{v}(\tilde{\varrho}_v;\beta,r)$.
\par   The relative increase of risk of the semiparametric
IC~$\tilde{\varrho}_v$ over that of the robust IC~$\eta_v$
remains to be investigated numerically---even in one dimension
when $\beta\ne \beta(r)$.
A suboptimal~$\tilde{\varrho}_v$ may still be useful.
\subsection{Discrepancy for Contamination}             \label{ss.r.c}
Contamination bias is $\omega_c(\psi)=\wsup_P|\psi|$,
the $L_{\infty}$~norm.
The robust IC~$\eta_c$ which minimizes~$\MSE_c(\ldotp;\beta,r)$,
by HR~(1994; Theorem~5.5.6), is the Hampel--Krasker influence curve,
\begingroup \mathsurround0em\arraycolsep0em \begin{eqnarray}
\label{e.etac} & \Ds \eta_c =
  (A \Lambda-a)\hspace{6\ei}w  \weg, \hspace{2em} w  =
  \min\Bigl\{ 1,\:\frac{b}{|A \Lambda-a|} \hspace{6\ei}\Bigr\} & \\
\noalign{\vspace{\belowdisplayshortskip}\pagebreak[0]\mathsurround=\msu
\noindent   with a particular bound, namely, the solution~$b$
to the equation
\nopagebreak\vspace{\abovedisplayskip}}
\label{e.boundc} & \Ds   \beta \hspace{6\ei}r^2 \hspace{4\ei}b =
  \hEw {\bigl(|{A \Lambda-a}|-b \bigr)_+}       &
\end{eqnarray}\endgroup
which may be nonunique only if $ \beta \hspace{6\ei}r = 0 $
(in which case $\eta_c = \hat{\varrho} $).
\par     The semiparametric IC~$\tilde{\varrho}_c$, by Theorem~\ref{t.rhoc},
has coordinates
\begin{equation} \label{e.r.spICc}
   \tilde{\varrho}_j= A^{j,1}\, (\Lambda_1 + r) \land u_1
   + \cdots + A^{j,k}\, (\Lambda_k + r) \land u_k
\end{equation}
with upper clipping constants~$u_j$ defined by~(\ref{e.gc}),
and $(A^{j,i})^{-1}={\cal K}$ by~(\ref{e.Kballs}).
\par    Thus, in general, $\tilde{\varrho}_c$ is unbounded so that the risk
  $\MSE_c(\tilde{\varrho}_c;\beta,r)$ becomes infinite if only
  $\beta \hspace{6\ei}r > 0 $ (the only interesting case).
\par    The intuitive convex combinations, which have been used in robust
statistics prior to any other type of neighborhoods, have always turned out
very similar to total variation in robustness respects. It is therefore
surprising that the semiparametric recipe
        (\ref{e.cICballs}),~(\ref{e.Kballs}) may give reasonable results
for one model but not the other.
\par
In the simplest testing context (one parameter, one-sided),
however, the discrepancy will again disappear; confer
Remark~\ref{r.s.proj.G10} and Theorem~\ref{t.S.sad.c} below.
\section{Unresolved: Robust Adaptive Estimation}
\setcounter{equation}{0}\label{s.O}
In the general semiparametric model of Section~\ref{s.P}, given the
canonical influence curves~(\ref{e.cIC}), one~$\varrho_{\theta,\nu}$
for each parameter $\theta\in \Theta$, $\nu\in H_{\theta}$,
the construction problem is to obtain an estimator~$(S_n)$ that, for each
$\theta\in \nolinebreak\Theta$ and $\nu\in \nolinebreak H_{\theta}$,
is asymptotically linear at~$(\theta,\nu)$ with prescribed
IC~$\varrho_{\theta,\nu}$.
\par     \bigskip  \noindent    {\bf Infinitesimal Nonrobustness}
Such estimators are automatically nonrobust in the same asymptotic,
infinitesimal, setup in which their efficiency is obtained.
\par     For example, consider the model
  $dQ_{\theta,\nu}(x)=\nu(x- \theta)\,dx$ with location parameter
  $\theta\in\R$ and nuisance parameter~$\nu$ any symmetric Lebesgue
density of finite Fisher information of location
  ${\cal I}_{\nu}=\int \Lambda_{\nu}^2(x)\hskip9\ei\nu(x)\,dx$ with
  $\Lambda_{\nu}=- \dot \nu/\nu$; then
  $\Lambda_{\theta,\nu}(x)= \Lambda_{\nu}(x - \theta)$.
In this model, adaptivity $\Pi_{2;\theta,\nu}(\Lambda_{\theta,\nu})=0$
holds by reasons of symmetry.
Adaptive estimators have been constructed by \mbox{Beran}~(1974)
and Stone~(1975) which, at each $(\theta,\nu)$, achieve
expansion~(\ref{e.Sal}) with influence curve
   $ \varrho_{\theta,\nu}(x) = \hat{\varrho}_{\theta,\nu}(x) =
     {\cal I}_{\nu}^{-1}\Lambda_{\nu}(x- \theta) $,
that is, are asymptotically linear with~IC~$\hat{\varrho}_{\theta,\nu}$.
Hence, under~$Q_{\theta,\nu}^n$, these estimators achieve the most
concentrated limit law~${\cal N}(0,{\cal I}_{\nu}^{-1})$ in~(\ref{e.ufoN}),
as if~$\nu$ was known.
\par     The assumption of exact symmetry, however, is rather strict.
In practice, one would accept a distribution function as symmetric
if it only is in a small neighborhood of an exactly symmetric one.
Such nonparametric hypotheses of approximate symmetry have been
investigated by HR~(1981\,a; Section~3) and generalized
         by Kakiuchi and Kimura~(2000).
If~$Q_{\theta,\nu}$ is thus enlarged to a shrinking neighborhood
$U_*(\theta,\nu; r/\!\sqrt{n}\,)$, while still~$\theta$ has to be estimated,
the adaptive estimates $\sqrt{n}\,(S_n -\theta)$ are driven off from their
limit~${\cal N}(0,{\cal I}_{\nu}^{-1})$ by some bias up to
   $\pm r \, \omega_*(\hat{\varrho}_{\theta,\nu})$ which, for gross error
neighborhoods ($*=v,c$), may become infinite if only
   $\Lambda_{\nu}=-\dot\nu/\nu$ is unbounded.
This readily follows from the asymptotic linearity~(\ref{e.Sal}) and
the results in HR~(1994; Section~5.3).
\par     The automatic nonrobustness of efficient estimators, under
asymmetric gross errors, in particular answers the question raised by
Huber~(1981; \S~1.2, p~7). The extension to the general semiparametric
model with unbounded canonical influence curve~$\varrho_{\theta,\nu}$ is
obvious.
\par     \bigskip  \noindent    {\bf Other Robustness Aspects}
Not considered here are qualitative robustness and (positive) breakdown
point. Like Huber~(1996; Section~28), we conjecture them to be incompatible
with adaptiveness, that is, asymptotic efficiency, for the usual
semiparametric models.
Possibly related is the necessary nonuniform convergence of adaptive
estimators in the symmetric location case; confer Bickel's~(1981)
presentation of Klaassen's result. Similar results by
Pfanzagl and Wefelmeyer~(1982; Proposition~9.4.1, Corollary~9.4.5)
connect this nonuniform convergence more explicitly with the
discontinuity of Fisher information. 
On the contrary, it is easy to see (since the Lindeberg condition may
be verified uniformly) that Huber's~(1964) minimax location M-estimate
tends to its normal limits uniformly on the corresponding symmetric
contamination neighborhood.
Thus it seems that a robustification would entail
other desirable properties.
\par     \bigskip \noindent
          {\bf Adaptation of Optimally Robust Estimators}
In view of all this, it seems desirable to construct estimators not with
the canonical influence curves~$\varrho_{\theta,\nu}$ but the robust
influence curves~$\eta_{\theta,\nu}$ instead, sacrificing a few percent
efficiency under each~$Q_{\theta,\nu}$ to gain robustness against
deviations from~$Q_{\theta,\nu}$.
\par     A first step in this direction has been made by
Shen~(1995; Theorem~2) who derives a bounded influence curve
$\eta_{\theta,\nu}=\eta_c$ minimizing $\Ew|\psi|^2$ among all
influence curves $\psi\in \Psi$, as defined in~(\ref{e.IC}) for a
general semiparametric model, subject to $|\psi|\le \sup|\eta_c|$.
In some sense, the result may be viewed an extension of
HR~(1994; Theorem~5.5.1), from finite to infinite dimensional nuisance
tangent space~$\partial_2 {\cal Q}$ of a certain kind; namely,
an~$L_2$~space of functions, expectation zero, and measurable relative
to some sub~$\sigma$~algebra of~${\cal B}$. For a similar result
and proof, confer HR~(1994; Theorem~7.4.13).
\par  
The corresponding adaptive estimator construction, however,
has not been achieved yet. The construction by Shen~(1994; Theorem~2)
in the symmetric location case is again only a first step, since
the desired asymptotic linearity and influence curve are established
but not the required uniform behavior of the estimator over shrinking
full neighborhoods.
\section{Semiparametric \protect \boldmath $C(\alpha)$-Tests}   %
\setcounter{equation}{0}                                      \label{s.T}
The semiparametric approach may further be applied to the testing of
hypotheses about the parameter of interest. The optimal tests are
generalized $C(\alpha)$-tests, which are based on residual scores
after an orthogonal projection on the closed linear tangent space
for the nuisance parameter.
In connection with the robust tangent balls, the nonlinear projection
on these balls will be employed instead, and yields sensibly bounded
modifications of the test statistics of the classical, asymptotically
maximin, multiparameter tests.
\subsection{\ran \protect \boldmath
            $C(\alpha)$-Tests For Tangent Spaces}    \label{ss:Tl}
We invoke the general setup of Section~\ref{s.P}:
The semiparametric probability model~${\cal Q}=\{\, Q_{\theta,\nu} \mid
                  \theta\in \Theta\weg,\:\nu\in H_{\theta}\,\}$
with main parameter~$\theta$, nuisance parameter~$\nu$,
the fixed parameter value~$(\theta_0,\nu_0)$ and corresponding element
      $ Q=\nolinebreak Q_{\theta_0,\nu_0} $, the scores
function~$\Lambda\in L_2^k=L_2^k(Q)$ of~${\cal Q}$ for~$\theta$ and the
differentiablity~(\ref{e.L2diff}) of~${\cal Q}$ at~$(\theta_0,\nu_0)$,
the orthogonal projection
  $\Pi_2\colon L_2^k\to (\clin\partial_2 {\cal Q})^k$,
and the Fisher information ${\cal J}=\Cov \bar{\Lambda}$ of~${\cal Q}$
for~$\theta$ at~$(\theta_0,\nu_0)$, 
where~$\bar{\Lambda}$ denotes the residual scores
\begin{equation}
         \bar{\Lambda}=\Lambda- \Pi_2(\Lambda)
\end{equation}
Given some numbers $ -\infty<z_1<z_2<\infty$ and $0\le z_3<z_4<\infty$,
local asymptotic one- and multisided hypotheses about the difference between
the true~$\theta$ and its reference value~$\theta_0$ are defined by
\begingroup \mathsurround0em\arraycolsep0em \begin{eqnarray}
   H': \hspace{24\ei}
   e' {\cal J}^{1/2} a \le z_1 \hspace{2em} &{\rm vs.}&  \hspace{2em}
   K': \hspace{24\ei} e' {\cal J}^{1/2} a \ge z_2 \\
   H'': \hspace{1.25em}
   a' {\cal J} a \le z_3^2 \hspace{2em} &{\rm vs.}& \hspace{2em}
   K'': \hspace{1.25em} a' {\cal J} a \ge z_4^2
\end{eqnarray}\endgroup
where $e\in\R^k$, $|e|=1$, is some fixed unit vector, and~${\cal J}^{1/2}=A$
any $k \times k$ root of~${\cal J}$ such that $AA'={\cal J}$.
\par
The hypotheses concern the sequence of laws~$Q_n\in {\cal Q}$
of the~$n$ i.i.d.\ observations $x_1,\ldots,x_n \sim Q_n$.
It is assumed that, for any $a\in\R^k$ and $g\in \partial_2 {\cal Q}$,
eventually 
\begin{equation}
    Q_n= Q_n(a,g)=Q_{ \theta_0+s_na,\,\nu_{s_n}^{\hspace{6\ei}g} }
\hspace{-1.5em}\samepage\end{equation}
where $ s_n=1/\!\sqrt{n}\,$ and $t \mapsto \nu_{t}^g\in H_{\theta_0+ta}$
is some path with tangent~$g$ in~(\ref{e.L2diff}).
\par
We employ asymptotic tests~$\delta=(\delta_n)$, that is, sequences of
tests~$\delta_n$ at sample size~$n$. Their error probabilities will be
evaluated under the $n$~fold product measures~$Q_n^n$, 
asymptotically, as $n$ tends to infinity.
\par
For~$\alpha\in (0,1)$, let~$u_{\alpha}$ denote the upper $\alpha$~point
of the standard normal distribution~$\Phi$, such that
  $\Phi(-u_{\alpha})=\alpha$. By~$\chi^2(k,z^2)$ denote the
  $\chi^2$~distribution with $k$~degress of freedom and noncentrality~$z^2$,
respectively a random variable having this distribution,
and by~$c_{\alpha}(k,z^2)$ its upper $\alpha$~point.
\pagebreak[2]\begin{Thm}\sl \label{t.t.Ca}
Let~$\delta=(\delta_n)$ be any sequence of tests.
\begin{ABC} \item    \label{i.t.Ca.1}
Then, in the one-sided case,
\begingroup \mathsurround0em\arraycolsep0em \begin{eqnarray}
\label{e.t.lev1}& \Ds   \sup_{H'} \limsup_{n\to \infty}
   \int \delta_n \,dQ_n^n(a,g) \le \alpha & \\
\noalign{\noindent implies\nopagebreak}
\label{e.t.pow1max} & \Ds \inf_{K'} \limsup_{n\to \infty}
   \int \delta_n \,dQ_n^n(a,g)
   \le \Phi \bigl(-u_{\alpha}+(z_2-z_1)\bigr)
\end{eqnarray}\endgroup
\item    \label{i.t.Ca.m}
In the multisided case,
\begingroup \mathsurround0em\arraycolsep0em \begin{eqnarray}
\label{e.t.levm} & \Ds \sup_{H''} \limsup_{n\to \infty}
   \int \delta_n \,dQ_n^n(a,g)    \le \alpha & \\
\noalign{\noindent implies\nopagebreak}
\label{e.t.powmmax} & \Ds \inf_{K''} \limsup_{n\to \infty}
  \int \delta_n \,dQ_n^n(a,g)
   \le \Pr \bigl(\, \chi^2(k,z_4^2) > c_{\alpha}(k,z_3^2)\bigr)
\end{eqnarray}\endgroup
\item    \label{i.t.Ca.att}
Bounds\/~{\rm (\ref{e.t.pow1max})} and\/~{\rm (\ref{e.t.powmmax}),}
with~$\limsup$ replaced by\/~{\rm $\liminf$,} are achieved by
the asymptotic tests
\begingroup \mathsurround0em\arraycolsep0em \begin{eqnarray}
   \delta'= (\delta'_n)\weg,\hspace{1.5em}   \delta'_n &{}={}&
   \hJc \bigl(e' {\cal J}^{-1/2} Z_n >  u_{\alpha} + z_1 \bigr) \\
   \delta''= (\delta''_n)\weg,\hspace{1.5em} \delta''_n &{}={}&
   \hJc \bigl( Z_n' \hspace{6\ei}{\cal J}^{-1} Z_n >
                           c_{\alpha}(k,z_3^2)    \bigr)
\end{eqnarray}\endgroup
respectively,
where  $ Z_n= 1/\!\sqrt{n}\,\sum_{i=1}^n \bar{\Lambda}(x_i) $.
\end{ABC} \end{Thm}
\begin{Bew}\enskip
The differentiability~(\ref{e.L2diff}),
for every $a\in\R^k$, $ g\in \partial_2 {\cal Q}$,
entails the following loglikelihood expansion,
\begin{equation} \label{e.s.llh}
  \log \frac{dQ_n^n(a,g)}{dQ^n} =
  \frac{1}{\sqrt{n}\,}\sum_{i=1}^n (a' \Lambda+g)(x_i)
  - \frac{1}{2} {\Vert a' \Lambda+g\Vert}^2 + \Lo_{Q^n}(n^0)
\end{equation}
Thus, given~$a\in\R^k$, the Fisher information
  $\Vert a'\Lambda+g\Vert^2$ at~$t=0$ of the one parameter family
  ${\cal Q}(a,g)=\{Q_{ \theta_0+ta,\,\nu_t^g}\}$ is minimized with
respect to $g\in \partial_2 {\cal Q}$ by
    $ g_a=- \pi_2(a' \Lambda)=-a' \hspace{9\ei} \Pi_2(\Lambda)$.
Therefore, associating with each $a\in\R^k$ any path~$\nu_t^a=\nu_t^{g_a}$,
the sequence of $k$~parameter submodels
    ${\cal Q}_n = \{\,Q_{n,a} \mid a\in\R^k\,\} $
consisting of the elements $Q_{n,a}=\nolinebreak Q_n(a,g_a)$,
will turn out least favorable.
\par     In fact, as $a' \Lambda+g_a=a' \bar{\Lambda}$ and
  $\Cov \bar{\Lambda}= {\cal J}$, expansion~(\ref{e.s.llh}) specializes to
\begin{equation}
    \log \frac{dQ_{n,a}^n}{dQ^n} = 
    \frac{a' }{\sqrt{n}\,}\sum_{i=1}^n \bar{\Lambda}(x_i)
  - \frac{1}{2} a' {\cal J}a + \Lo_{Q^n}(n^0)
\end{equation}
Because of this asymptotic normality, of the sequence of 
product models~${\cal Q}_n^n$, Theorems~3.4.6,~3.4.11 of HR~(1994) are
in force and, subject to (\ref{e.t.lev1}) and~(\ref{e.t.levm}),
respectively, furnish the power bounds (\ref{e.t.pow1max})
and~(\ref{e.t.powmmax}), as well as the asymptotically most powerful
level~$\alpha$ tests $\delta'$ and~$\delta''$, for the sequence of
submodels. 
\par
But, for arbitrary tangents~$g\in \partial_2 {\cal Q}$, 
(\ref{e.s.llh}) implies the following asymptotic normality of~$Z_n$
under~$Q_n^n(a,g)$,
\begingroup \mathsurround0em\arraycolsep0em \begin{eqnarray} &\Ds
   Z_n(Q_n^n(a,g)) \gwto {\cal N} \bigl(\hspace{6\ei}
   \Ew \bar{\Lambda}(a' \Lambda+g), {\cal J}\hspace{6\ei} \bigr)  &\\
\noalign{\noindent where\nopagebreak} &\Ds
  \Ew \bar{\Lambda}(a' \Lambda+g) = \Ew \bar{\Lambda}\Lambda'a
  = {\cal J}a \hspace{-36\ei} &
\end{eqnarray}\endgroup
since~$\bar{\Lambda}$ is orthogonal to~$\Pi_2(\Lambda)$ and~$g$.
Hence, the asymptotic error probabilities of the tests~$\delta'$
and~$\delta''$ do not depend on~$g\in \partial_2 {\cal Q}$.
\qed\end{Bew}
\begin{Rem}\rm
The orthogonality of~$\bar{\Lambda}$ on~$\partial_2 {\cal Q}$ may be used
a second time to construct test statistics that do not require knowledge
of~$\nu_0$. 
\par     In the finite dimensional case, confer Remark~\ref{r.p.fidi},
upon a regularization of the likelihoods,
(total) scores function,   
and Fisher information, estimates of~$\nu$ which are $\sqrt{n}\,$~consistent
and suitably discretized may be inserted for~$\nu_0$;
confer HR~(1994; Lemmas~6.4.1 and~6.4.4). 
Thus Neyman's $C(\alpha)$-tests are obtained, under no stronger conditions
than mean square differentiable root densities and identifiability
(of the ideal model). 
\par
The test statistics $ Z_n $~may also be replaced by an estimator $S=(S_n)$
of~$\theta$ which is asymptotically linear, in the sense of~(\ref{e.Sal}),
at each~$Q_{\theta,\nu}$, with canonical influence
curve~$\varrho_{\theta,\nu} = {\cal J}^{-1}_{\theta,\nu}
                 \bar{\Lambda}_{\theta,\nu} $;
confer HR~(1994; Theorem~6.4.8). This leads to
Wald's estimator tests $\lambda'$ and~$\lambda''$,
\begingroup \mathsurround0em\arraycolsep0em \begin{eqnarray}
\hspace{-1.5em}   \lambda'=(\lambda'_n)\weg,\hspace{1.5em}
   \lambda'_n &{}={}&  \hJc \bigl(
   e' {\cal J}^{1/2} \hspace{6\ei}\sqrt{n}\,(S_n- \theta_0)
                      >  u_{\alpha} + z_1 \bigr) \\
\hspace{-1.5em}   \lambda''=(\lambda''_n)\weg,\hspace{1.5em}
   \lambda''_n  &{}={}&  \hJc \bigl(
    n \,(S_n- \theta_0)' \hspace{6\ei} {\cal J} \,(S_n- \theta_0)
    >  c_{\alpha}(k,z_3^2)    \bigr)
\end{eqnarray}\endgroup
Like~$\delta'$ and~$\delta''$, also the test sequences $\lambda'$
and~$\lambda''$ achieve maxmin asymptotic power subject to level~$\alpha$
for~$H'$ vs.~$K'$, respectively for~$H''$ vs.~$K''$.
\par     In the infinite dimensional case, the estimation
of~$ \bar{\Lambda}_{\theta_0,\nu_0}$ and~${\cal J}_{\theta_0,\nu_0}$
(with $\theta_0$ known, $\nu_0$ unknown), and the construction of an
asymptotically linear estimator with canonical influence curve 
  $ \varrho_{\theta,\nu} = {\cal J}^{-1}_{\theta,\nu}
    \bar{\Lambda}_{\theta,\nu} $ at~$Q_{\theta,\nu}$ (at least for
  $\theta=\theta_0$ and every $\nu\in H_{\theta_0}$) is more difficult.
The methods of Klaassen~(1987) and the references mentioned therein may
prove useful.\qed \end{Rem}
\subsection{\ran \protect \boldmath
            $C(\alpha)$-Tests For Tangent Balls}     \label{ss:Tb}
As in Section~\ref{s.B}, we start
from~${\cal P}= \{\,P_{\theta}\mid \theta\in \Theta\,\}$,
an ideal, smooth $k$~parametric model without nuisance parameter.
\subsubsection{Parametric Tests}
Theorem~\ref{t.t.Ca} first specializes with $\partial_2 {\cal P}=\{0\}$.
Thus, the classical test sequences 
\begingroup \mathsurround0em\arraycolsep0em \begin{eqnarray}
\label{e.t.pCa1}   \hspace{-1.5em} \hat{\delta}'=(\hat{\delta}'_n)\weg,
\hspace{1.5em} \hat{\delta}'_n &{}={}&
   \hJc \bigl( e' \hspace{6\ei}{\cal I}^{-1/2}
   \hat{Z}_n >  u_{\alpha} + z_1 \bigr)  \\
\noalign{\noindent and\nopagebreak}\label{e.t.pCam}
\hspace{-1.5em}\hat{\delta}''=(\hat{\delta}''_n)\weg, \hspace{1.5em}
 \hat{\delta}''_n &{}={}&  \hJc \bigl(
   \hat{Z}_n' \hspace{9\ei}{\cal I}^{-1} \hat{Z}_n >
                           c_{\alpha}(k,z_3^2)    \bigr)
\end{eqnarray}\endgroup
based on $ \hat{Z}_n = 1/\!\sqrt{n}\,\sum \Lambda(x_i) $, as well as
the test sequences $\hat{\lambda}'$ and $\hat{\lambda}''$ employing an
asymptotically linear estimator $\hat{S}=(\hat{S}_n)$ with influence
curve~$\hat{\varrho}={\cal I}^{-1}\Lambda$ at~$P=P_{\theta_0}$,
\begingroup \mathsurround0em\arraycolsep0em \begin{eqnarray}
\noalign{\vspace{-\abovedisplayskip}\vspace{\abovedisplayshortskip}}
\label{e.t.pWa1}  \hat{\lambda}'_n &{}={}&
  \hJc \bigl( e' \hspace{6\ei}
  {\cal I}^{1/2} \hspace{6\ei}\sqrt{n}\,(\hat{S}_n- \theta_0)
      >  u_{\alpha} + z_1 \bigr) \\
\label{e.t.pWam}  \hat{\lambda}''_n &{}={}&  \hJc \bigl(
  n \,(\hat{S}_n- \theta_0)' \hspace{6\ei} {\cal I} \,(\hat{S}_n- \theta_0)
  >  c_{\alpha}(k,z_3^2)    \bigr)  \hspace{-1.5em}
\end{eqnarray}\endgroup
achieve maxmin asymptotic power subject to level~$\alpha$,
for the following parametric local asymptotic one- and multisided hypotheses
about $\theta-\theta_0$, respectively,
\begingroup \mathsurround0em\arraycolsep0em \begin{eqnarray}
\label{e.t.pHK1}   \widehat{H}': \hspace{24\ei}
   e' \hspace{6\ei}{\cal I}^{1/2} a \le z_1 \hspace{2em} &{\rm vs.}&
  \hspace{2em} \widehat{K}': \hspace{24\ei}
  e' \hspace{6\ei}{\cal I}^{1/2} a \ge z_2 \\
\label{e.t.pHKm}   \widehat{H}'': \hspace{1.25em} a' \hspace{6\ei}
  {\cal I} a \le z_3^2 \hspace{2em} &{\rm vs.}& \hspace{2em}
  \widehat{K}'': \hspace{1.25em} a' \hspace{6\ei}{\cal I} a \ge z_4^2
\end{eqnarray}\endgroup
\subsubsection{Semiparametric Projection}
Now enlarge the parametric measures~$P_{\theta}$ to
neighborhoods~$ U(\theta;r_0) $ under the null hypothesis,
respectively~$ U(\theta;r_1) $ under the alternative. Thus, robust
local asymptotic one- and multisided hypotheses $\widetilde{H}'$
vs.~$\widetilde{K}'$, and $\widetilde{H}''$ vs.~$\widetilde{K}''$
about $\theta-\theta_0$ are obtained.     
These concern the laws~$Q_n\in U(\theta_0+s_na;s_nr_{0/1})$ at
sample size~$n$, where $s_n=1/\!\sqrt{n}\,$, and $a\in\R^k$ is subject
to the conditions of the corresponding parametric hypotheses
  $\widehat{H}'$, $\widehat{K}'$, $\widehat{H}''$,~$\widehat{K}''$,
respectively. 
\par
By this enlargement, size and power of the tests $\hat{\delta}'$
and~$\hat{\delta}''$ will be affected without control.
Conceptually, a robustification is appealing that interpretes
model deviations as nuisance parameter.
Then, to the neighborhood model~${\cal Q}$ of semiparametric
form~(\ref{e.b.nbdQ}), Theorem~\ref{t.t.Ca} may again be applied,
and leads to the semiparametric recipe: From~$\Lambda$ subtract the
component~$\Pi_2(\Lambda)$ explained by the nuisance parameter,
and exchange the test statistics $ {\cal I}^{-1/2}\hat{Z}_n$
based on~$\Lambda$ for the test statistics $ {\cal J}^{-1/2}Z_n$
based on $ \bar{\Lambda}=\Lambda-\Pi_2(\Lambda)$.
\begin{Rem}\rm  In the context of testing, contrary to estimation, there
is no Fisher consistency requirement, that is,
  $ \Ew \psi \Lambda' =\EM_k $ in~(\ref{e.IC}) and the corresponding
  standardization by~$ {\cal J}^{-1}$ in~(\ref{e.cIC}).
The present standardization of~$\bar{\Lambda}$ by~${\cal J}^{-1/2}$ shall
achieve unit covariance of the limit normals to obtain invariance under
the orthogonal group, which is needed in the proof of the maxmin testing
result.\qed \end{Rem}
Because, in the case of Hellinger, total variation, and contamination
neighborhoods, the tangent sets $ \partial_2 {\cal Q}_* $, $*=h,v,c$,
determined by Proposition~\ref{p.balls} 
achieve $\clin \partial_2 {\cal Q}_*=L_2\cap \{\Ew=0\}$,
we replace, as we did in Section~\ref{s.C},
  $\pi_2$ and~$\Pi_2$ by the nonlinear projection
  $\widetilde{\pi}_2 \colon L_2\to \nolinebreak {\cal G}_* $
  on~$\partial_2 {\cal Q}_*={\cal G}_*$, respectively by
  $ \widetilde{\Pi}_2 ={(\widetilde{\pi}_2,\ldots,\widetilde{\pi}_2)}'
    \colon L_2^k\to {\cal G}_*^k$ (acting coordinatewise).  %
\par
Actually, the situation is more complex for testing than for estimation
in Section~\ref{s.C}, since now two neighborhoods (null hypothesis,
alternative) are involved.
This will be clarified in Remark~\ref{r.s.proj.G10} below.
 \pagebreak[1] \par
  We first put $r=r_0+r_1$ and naively project on~${\cal G}_*$
(of this radius~$r$). Thus, let
\begingroup \mathsurround0em\arraycolsep0em \begin{eqnarray}
\noalign{\vspace{-\abovedisplayskip}\vspace{\abovedisplayshortskip}}
\label{e.t.rob.L}
      \tilde{\Lambda} &{}={}& \Lambda- \widetilde{\Pi}_2(\Lambda) \\
\noalign{\noindent and suppose that\nopagebreak}
\label{e.t.rob.J}
      \tilde{{\cal J}} &{}={}& \Cov \tilde{\Lambda} > 0
\end{eqnarray}\endgroup
Then, based on $ \tilde{Z}_n = 1/\!\sqrt{n}\,\sum \tilde{\Lambda}(x_i) $,
the semiparametric approach leads to the scores statistics,
\begin{equation} \label{e.t.sp.stat.Ca}
    e' \tilde{{\cal J}}^{-1/2} \tilde{Z}_n \weg,\qquad
    \tilde{Z}_n' \hspace{6\ei} \tilde{{\cal J}}^{-1} \tilde{Z}_n
\end{equation}
for testing the robust one- and multisided hypotheses $\widetilde{H}'$
vs.~$\widetilde{K}'$ and $\widetilde{H}''$ vs.~$\widetilde{K}''$,
respectively. The corresponding semiparametric estimator tests would
employ the statistics
\begin{equation}  \label{e.t.sp.stat.Wa}
   e' \tilde{{\cal J}}^{-1/2}\hspace{6\ei}{\cal K}\,
   \sqrt{n}\,(\tilde{S}_n- \theta_0)\weg,\qquad
   n \,(\tilde{S}_n- \theta_0)' \hspace{6\ei} {\cal K}'
   \tilde{{\cal J}}^{-1} \hspace{6\ei}{\cal K}\,(\tilde{S}_n- \theta_0)
\end{equation}
based on an asymptoticaly linear estimator $\tilde{S}=(\tilde{S}_n)$
with semiparametric influence curve
  $\tilde{\varrho}={\cal K}^{-1}\tilde{\Lambda}$,  
provided ${\cal K}= \Ew \tilde{\Lambda}\Lambda'$ is regular; confer
  (\ref{e.cICballs}),~(\ref{e.Kballs}).
\par     The semiparametric asymptotic tests thus obtained are denoted by
  $ \tilde{\delta}'$, $ \tilde{\delta}''$, and~$ \tilde{\lambda}'$,
$ \tilde{\lambda}''$, respectively. The suitable choice of the critical
values for their test statistics, however, must be left open. 
\subsubsection{`Robust' Test Statistics}
\paragraph{Hellinger Model}
By Theorem~\ref{t.rhoh}, under condition~(\ref{e.rh}):
   $ 8 \,r^2 < \min\, {\cal I}_{j,j} $,
we have $ \tilde{\Lambda}=D \Lambda $ with regular matrix
        $D=\diag(1- \gamma_j)$, where
        $ \gamma_j^2 \hspace{9\ei}{\cal I}_{j,j} = 8 \,r^2 $.
\par    It follows that $ \tilde{{\cal J}}=D\,{\cal I}D$,
   $ \tilde{{\cal J}}^{1/2}= D\,{\cal I}^{1/2}$, and so
   $ \tilde{{\cal J}}^{-1/2} \tilde{\Lambda}= {\cal I}^{-1/2}\Lambda$.
Moreover, ${\cal K}= D \,{\cal I}$, hence
    $ {\cal K}' \tilde{{\cal J}}^{-1} \hspace{6\ei}{\cal K} = {\cal I}$,
and $ \tilde{\varrho}=\hat{\varrho}= {\cal I}^{-1}\Lambda$
by Theorem~\ref{p.etHro}.
\par     Therefore, the semiparametric test statistics
         (\ref{e.t.sp.stat.Ca}),~(\ref{e.t.sp.stat.Wa}) agree with the
parametric test statistics in~(\ref{e.t.pCa1})--(\ref{e.t.pWam}).
         The result matches Theorem~\ref{p.etHro}.
\paragraph{Total Variation Model}
Under condition~(\ref{e.rv}): $ 2 \,r < \min \,\Ew |\Lambda_j| $,
Theorem~\ref{t.rhov} furnishes~$ \tilde{\Lambda} $ with
coordinates~$ \tilde{\Lambda}_j = v'_j \lor \Lambda_j \land v''_j $
and clipping constants determined by~(\ref{e.gv}).
Thus the coordinates of~$ \tilde{{\cal J}}^{-1/2} \tilde{\Lambda} $ are
linear combinations of~$v'_j \lor \Lambda_j \land v''_j$, hence are bounded.
\par     Boundedness of the semiparametric test statistics and
influence curve~$\tilde{\varrho}_v$, confer~(\ref{e.r.spICv}),
ensures a minimal robustness of the corresponding semiparametric tests
  $ \tilde{\delta}_v'$, $ \tilde{\lambda}'_v$, for $\widetilde{H}'_v$
vs.~$\widetilde{K}'_v$, and $ \tilde{\delta}_v''$, $ \tilde{\lambda}_v''$
for $\widetilde{H}''_v$ vs.~$\widetilde{K}''_v$.
\paragraph{Contamination Model}
Under condition~(\ref{e.rc}):
  $ r < - \max \,\winf_{P} \Lambda_j $, Theorem~\ref{t.rhoc} supplies
  $ \tilde{\Lambda}_j = (\Lambda_j + r) \land u_j $, whose upper clipping
constant~$u_j$ is defined by~(\ref{e.gc}). Thus the coordinates
of~$ \tilde{{\cal J}}^{-1/2} \tilde{\Lambda} $, certain linear
combinations of~$ (\Lambda_j + r) \land u_j $, may be unbounded.
\par     Unboundedness of the semiparametric test statistics and
influence curve~$\tilde{\varrho}_c$, confer~(\ref{e.r.spICc}),
entails maximum asymptotic error probabilities~$100\%$ of the
corresponding tests for the robust hypotheses;
as with estimation in Subsection~\ref{ss.r.c}.
\par
However, Remark~\ref{r.s.proj.G10} tells us that, instead on
  $ {\cal G}_c = r \hspace{6\ei}G_c $, we must actually project
on the set $ r_0 \hspace{6\ei}G_c - r_1 \hspace{6\ei}G_c $ (which
makes no difference in the Hellinger and total variation models.)
The correct~$\tilde{\Lambda}$ and~$\tilde{\varrho}_c$, therefore,
are determined by Theorem~\ref{t.S.sad.c}, and turn out bounded
towards both sides.
\par     Boundedness of the semiparametric test statistics and
influence curve~$\tilde{\varrho}_c$, now essentially of
form~(\ref{e.r.spICv}), again ensures some minimal robustness of the
corresponding tests
    $ \tilde{\delta}_c'$, $ \tilde{\lambda}'_c$ for $\widetilde{H}_c'$
vs.~$\widetilde{K}_c'$, and $ \tilde{\delta}_c''$, $ \tilde{\lambda}_c''$
for $\widetilde{H}_c''$ vs.~$\widetilde{K}_c''$.
\paragraph{Multiparameter, Multisided Case}
In this general case, an exact evaluation of the asymptotic maximum size
over~$\widetilde{H}''$ and minimum power over $\widetilde{K}''$ of the
derived semiparametric tests, and other tests based on quadratic forms
in sums or in asymptotically linear estimators, is rather complicated;
confer HR~(1994; \S~5.4, pp 192--194), especially equation~(54) there.
Optimization problems arise for the maximum eigenvalue of the information
standardized covariance subject to bounds on the self-standardized
sensitivity; see equation~(55), p~194. 
\subparagraph{Dimensional Advantage}
As these problems have not been solved yet, no optimally robust test is
distinguished, in comparison to which the semiparametric tests might be
judged.
\par     It certainly is an advantage of the semiparametric over the
maxmin approach to robust testing that it works in higher dimensions as
it works in one, and that it yields test statistics which seem reasonably,
if not optimally, robust.
\paragraph{One Parameter, One-Sided Case}
In the simplest case, a strong justification of the semiparametric approach
is possible. Section~\ref{s.S} will establish optimal robustness:
For the one parameter, one-sided, robust hypotheses $\widetilde{H}'$
vs.~$\widetilde{K}'$, the semiparametric test~$ \tilde{\delta}'$
(and~$ \tilde{\lambda}'$) is asymptotically maxmin.
\section{Saddle Points For Testing Convex Sets}
\setcounter{equation}{0}                                      \label{s.S}
Consider hypotheses which consist of local alternatives generated by any
two disjoint closed convex sets $G_0$ and~$G_1$ of tangents at some
probability~$P$. Picking the unique minimum norm element of $G_1-G_0$,
and the corresponding sequence of Neyman--Pearson tests, seems to fit
the semiparametric projection arguments---and furnishes a saddle point.
\par
The result applies to infinitesimal Hellinger, total variation, and
contamination neighborhoods around~$P$ and a local alternative of~$P$
with fixed tangent, respectively.
In the total variation and contamination cases, the maxmin asymptotic tests
thus obtained by projection agree with the robust asymptotic tests based
on the least favorable pairs in the sense of Huber and Strassen~(1973).
\subsection{Convex Sets Defining Local Alternatives}
Let~$P\in {\cal M}$ be some probability. Every tangent
    $\rho\in L_2\cap \{\Ew=0\}$ at~$P$ gives rise to a sequence of
local alternatives~$P_{n,\rho}$ of~$P$ such that, in the Hilbert space
of square root densities,
\begingroup \mathsurround0em\arraycolsep0em \begin{eqnarray}
\label{e:t:pathdiff}   \sqrt{dP_{n,\rho}}\, & {}={} &
  \bigl(1+ \Tfrac{1}{2}s_n \hspace{6\ei} \rho \bigr) \sqrt{dP}\,
  + \Lo(s_n) \qquad \mbox{as $\Ds n\to \infty$}\hspace{-2em}\\
\noalign{\vspace{\belowdisplayshortskip}
\mathsurround=\msu\noindent
    where $s_n=1/\!\sqrt{n}\,$.
    Constructions to achieve~(\ref{e:t:pathdiff})
    are\pagebreak[0]\vspace{\abovedisplayskip}} 
\label{E:I:pathdef}    \frac{dP_{n,\rho}}{dP} & {}={} &
   \cases{ \Bigl(\Tfrac{1}{2}s_n \hspace{6\ei} \rho + \sqrt{1-
           \Tfrac{1}{4}s_n^2 \hspace{6\ei}\Vert \rho\Vert^2 \hspace{6\ei}}
           \hspace{24\ei}\Bigr)^2\weg, & or simply\cr \hspace{1.5em}
           1+ s_n \hspace{6\ei} \rho & if~$\rho\in L_{\infty}$}
\end{eqnarray}\endgroup
Let $G_0,G_1\subset L_2\cap \{\Ew=0\}$ be any two disjoint sets of tangents.
The observations $x_1,\ldots,x_n$ at sample size~$n$ are assumed independent
identically distributed with distribution~$Q_n$.
For fixed $g=(g_0,g_1)\in G_0 \times G_1$, preliminary simple asymptotic
hypotheses concerning~$Q_n$ are that, eventually,
\begin{equation} \label{e.S.HKg}
     H_{g_0}: Q_n = P_{n,g_0} \qquad  K_{g_1}: Q_n = P_{n,g_1}
\end{equation}  
As in Section~\ref{s.T}, asymptotic tests~$\delta=(\delta_n)$, that is,
sequences of tests~$\delta_n$ at sample size~$n$, are employed, and
their error probabilities are evaluated under the $n$~fold product
measures~$Q_n^n$.
\par     Then the testing problem $H_{g_0}$ vs.~$K_{g_1}$
at level~$\alpha\in (0,1)$,
\begingroup \mathsurround0em\arraycolsep0em \begin{eqnarray}
  \liminf_{n\to \infty} \int \delta_n \, dP_{n,g_1}^n &{}={}& \max{!} \\
\noalign{\noindent \raisebox{0ex}[0ex][0ex]{subject to}\nopagebreak}
  \limsup_{n\to \infty} \int \delta_n \, dP_{n,g_0}^n &{}\le{}& \alpha
\end{eqnarray}\endgroup
has the solution $\delta_{g}=(\delta_{n,g})$,
\begin{equation}
   \delta_{n,g}=\lJc \biggl \lgroup \frac{1}{\sqrt{n}\,}\sum_{i=1}^n
   g_{10}(x_i) > \Vert g_{10}\Vert \hspace{6\ei} u_{\alpha}
   + \langle g_{10}|g_0\rangle  \biggr \rgroup
\end{equation}
where $g=(g_0,g_1)$, $g_{10}=g_1-g_0$, and~$u_{\alpha}$ denotes
the standard normal upper $\alpha$~point.
Under~$H_{g_0}$, $\delta_{g}$~achieves asymptotic size~$\alpha$ and
under~$K_{g_1}$, asymptotic power~$\Phi(-u_{\alpha}+\Vert g_{10}\Vert\,)$.
The tests~$\delta_{n,g}$ are unique up to terms tending to~$0$
in $P^n$~probability. All these statements follow from the
loglikelihood expansion~(\ref{e.S.llhg10}) below and
HR~(1994; Corollary~3.4.2\footnote{Note that $\sigma>0$
  must be assumed in part~(b).}).
\par     Put  $ H_{G_0} = \cup\{H_{g_0}\mid g_0\in G_0\} $ and
              $ K_{G_1} = \cup\{K_{g_1}\mid g_1\in G_1\} $.
\subsection{The Maxmin Test Result}
Then the maxmin testing problem $H_{G_0}$ vs.~$K_{G_1}$
at level $\alpha\in (0,1)$ is
\begingroup \mathsurround0em\arraycolsep0em \begin{eqnarray}
   \inf_{g_1\in G_1}\liminf_{n\to \infty}
   \int \delta_n \, dP_{n,g_1}^n  &{}={}& \max {!} \\
\noalign{\noindent \raisebox{0ex}[0ex][0ex]{subject to}\nopagebreak}
   \sup_{g_0\in G_0}\limsup_{n\to \infty}
   \int \delta_n \, dP_{n,g_0}^n  &{}\le{}& \alpha  \end{eqnarray}\endgroup
\paragraph{Convex Closed Tangent Sets}
The tangent sets
\begin{equation}
   G_0,G_1\subset L_2\cap \{\Ew=0\}\weg,\qquad G_0\cap G_1=\emptyset
\end{equation}
are each assumed convex and closed in~$L_2$. The set of differences
\begin{equation} \label{e.S.G10cl}
         G_{10}=G_1-G_0 
\end{equation}
which is again convex, may not be closed if $\dim L_2>2$, and therefore
explicitly assumed to be also closed. Then denote by~$q_{10}=q_1-q_0$
the unique minimum norm element of~$G_{10}$; as $G_0\cap G_1=\emptyset$,
we have $q_{10}\ne0$.
\begin{Thm}\sl \label{t.x}
The asymptotic testing problem $H_{G_0}$ vs.~$K_{G_1}$ at level~$\alpha$
has a saddle point at\/~{\rm $q=(q_0,q_1)$,}
and the maxmin asymptotic power achieved by~$\delta_{q}$
equals~$\Phi (-u_{\alpha}+\Vert q_{10}\Vert\,)$.
\end{Thm}
\begin{Bew}\enskip
For any tangent~$\rho$, the following loglikelihood expansion holds,
\begingroup \mathsurround0em\arraycolsep0em \begin{eqnarray}
\label{e.S.llhg}
   \log \frac{dP_{n,\rho}^n}{dP^n} &{}={}& s_n\Nsum_{i} \rho(x_i)
   -  \Tfrac{1}{2} \hspace{6\ei}{\Vert \rho\Vert}^2 + \Lo_{P^n}(n^0) \\
\noalign{\noindent \raisebox{0ex}[0ex][0ex]{Hence}\nopagebreak}
\label{e.S.llhg10}
    \log \frac{dP_{n,g_1}^n}{dP_{n,g_0}^n} &{}={}&
    s_n\Nsum_{i} g_{10}(x_i) + \mbox{const} + \Lo_{P^n}(n^0)
\end{eqnarray}\endgroup
by mutual contiguity, for every $g=(g_0,g_1)\in G_0 \times G_1$ and
$g_{10}=g_1-g_0$. Therefore, 
the test sequence~$\delta_{g}$ is indeed the optimum one at level~$\alpha$
for~$H_{g_0}$ vs.~$K_{g_1}$; confer HR~(1994; Corollary~3.4.2).
\par
Let us evaluate~$\delta_{q}$, for any $q=(q_0,q_1)\in G_0 \times G_1$ fixed,
under other tangents~$\rho\in G_0\cup G_1$. In view of (\ref{e.S.llhg}),
by LeCam's third lemma, confer HR~(1994; Corollary~2.2.6), the
sequence of test statistics $ s_n\sum_{i} q_{10}(x_i) $
are asymptotically normal under~$P^n_{n,\rho}$,
\begingroup \mathsurround0em\arraycolsep0em \begin{eqnarray}
   s_n\Nsum_{i} q_{10}(x_i) &{}\gwto{}&
   {\cal N}\bigl( \hspace{6\ei} \langle q_{10}|\rho \rangle,
   \Vert q_{10}\Vert^2 \hspace{9\ei}\bigr) \\
\noalign{\noindent \raisebox{0ex}[0ex][0ex]{hence}\nopagebreak}
\hspace{-1.5em} \label{e.S.q10/rho}
   \lim_{n\to \infty}\int \delta_{n,q}\, dP^n_{n,\rho}
   \hspace{-24\ei} &{}={}& \hspace{-24\ei}
   \Phi \Bigl( -u_{\alpha} + \frac{ \langle q_{10}| \rho-q_0
   \rangle}{\Vert q_{10}\Vert }\,\Bigr)
\end{eqnarray}\endgroup
Therefore, the asymptotic size under $g_0\in G_0$ becomes maximal
at $g_0=q_0$, and the asymptotic power under $g_1\in G_1$ becomes
minimal at $g_1= q_1$, iff
\begin{equation} \label{e.S.q10char}
    \langle q_{10}| q_{10}-g_{10} \rangle\le 0 \qquad
    \forall\,g_{10}\in G_{10}
\end{equation}
By Lemma~\ref{l.c.approx}, this characterizes the minimum norm
element~$q_{10}$ of~$G_{10}$.\qed\end{Bew}
While $q_{10}=q_1-q_0$ is unique, there may exist other least favorable
pairs of tangents $g=(g_0,g_1)$ in~$G_{0}\times G_{1}$ achieving the
same $g_{10}=g_1-g_0=q_{10}$ of minimum norm in~$G_{10}=G_1-G_0$.
But then $\delta_{g}=\delta_{q}$, by the following corollary. So the maxmin
asymptotic level~$\alpha$ test for~$H_{G_0}$ vs.~$K_{G_1}$ is unique.
\begin{Cor}\sl
Let $g=(g_0,g_1)$ and~$q=(q_0,q_1)$ be two least favorable tangent pairs
in\/~{\rm $G_{0}\times G_{1}$.} Then
\begin{equation} \label{e.S.q10q=q10g}
  \langle q_{10}|g_{0}\rangle = \langle q_{10}|q_{0}\rangle \weg, \qquad
  \langle q_{10}|g_{1}\rangle = \langle q_{10}|q_{1}\rangle
\end{equation}
\end{Cor}
\begin{Bew}\enskip
By the saddle point, $\delta_{q}$~achieves asymptotic size $\le \alpha$
under~$H_{g_0}$ and asymptotic power
  $ \ge \Phi(-u_{\alpha}+ \Vert q_{10}\Vert \,) =
        \Phi(-u_{\alpha}+ \Vert g_{10}\Vert \,) $ under~$K_{g_1}$.
However, strict inequalities cannot hold since~$\delta_{g}$ is optimal for
$H_{g_0}$ vs.~$K_{g_1}$. Inserting $\rho=g_{0},g_{1}$ in~(\ref{e.S.q10/rho})
and~(\ref{e.S.q10char}), (\ref{e.S.q10q=q10g})~follows. Hence,
in particular, $\delta_{g}=\delta_{q}$.
\qed\end{Bew}
\subsection{Robust Asymptotic Tests}                       \label{ss:t.x}
In the setup of Section~\ref{s.B}, with $P=P_{\theta_0}$,
the normed robust tangent balls~$G_*$ are
\begingroup \mathsurround0em\arraycolsep0em \begin{eqnarray}
\noalign{\vspace{-\abovedisplayskip}\vspace{\abovedisplayshortskip}}
\label{e.S.Gh}   G_h & {}={} & \bigl\{\,g\in L_2 \bigm| \Ew g=0\weg,\;
    \Ew g^2\le 8 \,\bigr\} \\
\label{e.S.Gv}   G_v & {}={} & \bigl\{\,g\in L_2 \bigm| \Ew g=0\weg,\;
    \Ew |g|\le 2 \,\bigr\} \\
\label{e.S.Gc}   G_c & {}={} & \bigl\{\,g\in L_2 \bigm| \Ew g=0\weg,\;
     g\ge-1\,\bigr\}
\end{eqnarray}\endgroup
Thus ${\cal G}_*=r \hspace{6\ei}G_*$ are the balls of radius~$r$
         introduced in (\ref{e.b.Gh})--(\ref{e.b.Gc}); $*=h,v,c$.
\par     We assume parameter dimension $k=1$. Invoke the scores
function~$\Lambda$ for the parameter~$\theta$ of the ideal model~${\cal P}$
at~$\theta_0$, and let numbers $r_0,r_1,\tau\in [\,0,\infty)$ be given.
Then Theorem~\ref{t.x} is going to be applied to the tangent sets
\begin{equation}
         G_{*,0}=r_0 \hspace{6\ei}G_*\weg,\qquad
         G_{*,1}=\tau \hspace{3\ei}\Lambda+ r_1 \hspace{6\ei}G_*
\end{equation}
\begin{Rem}\rm \label{r.s.proj.G10}
The minimum norm element~$q_{*,10}$ of $G_{*,10}= G_{*,1}- G_{*,0}$,
therefore, will be $ \tau \hspace{3\ei}\Lambda $ minus its projection
on the set of differences
  $ r_0 \hspace{6\ei}G_* - r_1 \hspace{6\ei}G_* $.\qed\end{Rem}
Abbreviate the corresponding hypotheses by $H_*=H_{G_{*,0}}$
and~$K_*=K_{G_{*,1}}$. As shown in the proof to Proposition~\ref{p.balls},
    $H_*$ and~$K_*$ represent the neighborhoods
    $U_*(\theta_0;s_n \hspace{3\ei}r_0)$ and
    $U_*(\theta_0+s_n \hspace{3\ei}\tau; s_n \hspace{3\ei}r_1)$
about~$P_{\theta_0}$ and~$P_{\theta_0+s_n\tau}$ of radii
$s_n \hspace{3\ei}r_0$ and~$s_n \hspace{3\ei}r_1$ respectively,
up to some~$\Lo(s_n)$ where $s_n=1/\!\sqrt{n}\,$. Put $r=r_0+r_1$.
\subsubsection{Maxmin Tests for Hellinger Balls}
\begin{Thm}\sl \label{t.S.sad.h} Let
\begingroup \mathsurround0em\arraycolsep0em \begin{eqnarray}
\noalign{\vspace{-\abovedisplayskip}\vspace{\abovedisplayshortskip}}
\label{e.S.rad.h}       8 \,r^2 &{}<{}&
   \tau^2\,{\cal I} \weg,\quad \mbox{where}\enskip
   {\cal I}= \Vert \Lambda\Vert^2  \\
\noalign{\vspace{\belowdisplayshortskip}\pagebreak[0]\mathsurround=\msu
\noindent        Then the least favorable tangent pair
           $q_h=(q_{h,0},q_{h,1})$ in~$G_{h,0}\times G_{h,1}$
is unique,\nopagebreak\vspace{\abovedisplayskip}}
      q_{h,0}=r_0 \hspace{6\ei}\gamma \hspace{3\ei}\Lambda\weg, \quad
      q_{h,1}   &{}={} &\tau \hspace{3\ei}\Lambda -
      r_1 \hspace{6\ei}\gamma \hspace{3\ei}\Lambda\weg, \quad
      \mbox{where}\enskip \gamma= \sqrt{8}\:\Vert \Lambda\Vert^{-1}
\end{eqnarray}\endgroup
The maxmin asymptotic level~$\alpha$ test~$\delta_{q_h}=(\delta_{n,q_h})$
for~$H_h$ vs.~$K_h$ is given by
\begin{equation}
   \delta_{n,q_h} =
   \lJc \biggl \lgroup \frac{1}{\sqrt{n \hspace{9\ei}{\cal I}}\,}
   \sum_{i=1}^n \Lambda(x_i) >  u_{\alpha} + \sqrt{8}\,r_0 \biggr \rgroup
\end{equation}
and achieves maxmin asymptotic power
  $ \Phi \bigl(-u_{\alpha}+ \tau \hspace{6\ei}\Vert \Lambda\Vert
    - \sqrt{8}\,r\,\bigr) $.
\end{Thm}
\begin{Bew}\enskip Since $G_h$ is symmetric convex, $ G_{10} =
    \tau \hspace{3\ei}\Lambda + r_1 \hspace{6\ei}G_h - r_0 \hspace{6\ei}G_h
  = \tau \hspace{3\ei}\Lambda - r \hspace{6\ei}G_h $,
and the minimum norm element~$q_{10}$ is supplied by
  $ q_0= r_0 \hspace{6\ei}\tilde{g} $,
  $ q_1= \tau \hspace{3\ei}\Lambda - r_1 \hspace{6\ei}\tilde{g} $,
where~$\tilde{g}\in G_h$ is the unique minimizer of
  $ \Vert \tau \hspace{3\ei}\Lambda-r \hspace{6\ei}g\Vert $
among all $g\in G_h$.
\par
The projection of~$\tau \hspace{3\ei}\Lambda$ on~$r \hspace{6\ei}G_h$ is
determined by Theorem~\ref{t.rhoh} and its proof, with~$\Lambda$ replaced
by~$\tau \hspace{3\ei}\Lambda$. Then condition~(\ref{e.S.rad.h}) coincides
with condition~(\ref{e.rh}), and is equivalent to
  $r \hspace{6\ei}\tilde{g}\ne \tau \hspace{3\ei}\Lambda$, that is,
  $G_{h,0}\cap G_{h,1}=\emptyset$.
Thus $\tilde{g}=\gamma \hspace{3\ei}\Lambda$.
\par     \nopagebreak
With $q_{10}=(\tau- \gamma \hspace{6\ei}r)\hspace{6\ei}\Lambda$ and
     $ \langle q_{10}|q_0\rangle = \Vert q_{10}\Vert \hspace{8\ei}
       \sqrt{8}\,r_0$, Theorem~\ref{t.x} applies.
\par
The pair $(q_0,q_1)$ is unique: $q_0=r_0 \hspace{6\ei}g_0$
and $q_1=\tau \hspace{3\ei}\Lambda + r_1 \hspace{6\ei}g_1$ for arbitrary
elements $ g_0,g_1\in G_h$ entails that
  $ r \hspace{6\ei}\tilde{g} = r_0 \hspace{6\ei}g_0 - r_1 \hspace{6\ei}g_1 $,
and then $ g_0=-g_1=\tilde{g} $ because
  $ \Vert g_0\Vert, \Vert g_1 \Vert \le \sqrt{8}\, = \Vert \tilde{g}\Vert $
and the norm is strictly convex. \qed \end{Bew}
Thus the least favorable tangents---multiples of~$\Lambda$---generate
local alternatives within the parametric model~${\cal P}$, and
the asymptotic maxmin test~$\delta_{q_h}$ agrees with the asymptotic most
powerful test for $(P^n_{\theta_0})$ vs.~$(P^n_{\theta_0+s_n \tau})$
at the smaller level $ \Phi(-u_{\alpha}-\sqrt{8}\,r_0) $. 
The result compares with Theorem~\ref{p.etHro} and Remark~\ref{r.r.Hopt}.
\paragraph{Least Favorable Pairs of Probabilities}
For Hellinger balls, least favorable pairs of probabilities in the sense
of Huber and Strassen~(1973) do not exist; confer Birg\'e~(1980).
\par
For total variation and contamination neighborhoods, such Huber--Strassen
pairs exist.
While the least favorable pairs are not unique, their likelihood and its
distribution under each of the two probabilities of least favorable pairs
is unique; confer HR~(1977).
The Neyman--Pearson tests based on the likelihoods of the product measures
of least favorable probability pairs furnish finite sample size,
hence also asymptotic, maxmin tests.
\par
The robust asymptotic tests derived from Huber--Strassen pairs have been
evaluated by Huber--Carol~(1970), HR~(1978), Wang~(1981), and Quang~(1985).
\subsubsection{Maxmin Tests for Total Variation Balls}
\begin{Thm}\sl \label{t.S.sad.v} Let
\begin{equation} \label{e.S.rad.v}
      2 \,r < \tau \hspace{6\ei}\Ew|\Lambda|
\end{equation}
\begin{ABC}\item    \label{i.S.maxminT.v}
Then a least favorable tangent pair $q_v=(q_{v,0},q_{v,1})$
     in~$G_{v,0}\times G_{v,1}$ is given by
\begingroup \mathsurround0em\arraycolsep0em \begin{eqnarray}
\noalign{\vspace{-\abovedisplayskip}\vspace{\abovedisplayshortskip}}
\label{e.S.qv} q_{v,0} &{}={}& r_0 \hspace{6\ei}\tilde{g}_v\weg,\qquad
    q_{v,1}= \tau \hspace{3\ei}\Lambda -r_1 \hspace{6\ei}\tilde{g}_v \\
\noalign{\noindent where\nopagebreak} \label{e.S.clip.Lambda.v}
    r\hspace{6\ei}\tilde{g}_v &{}={}&
    \tau \hspace{6\ei}(\Lambda - v'' \hspace{6\ei})_+
    - \tau \hspace{6\ei}(v'- \Lambda)_+
\end{eqnarray}\endgroup
with clipping constants
   $v'=v'(r/\tau)<0<v''(r/\tau)=v''$ determined by
\begin{equation} \label{e.S.clip.c.v}
   \tau \hspace{3\ei}\hEw ( v' - \Lambda)_{+} = r
   = \tau \hspace{3\ei}\hEw ( \Lambda - v'' \hspace{6\ei})_{+}
\end{equation}
Setting\/ {\rm $ \Lambda^{(v)} = v'\lor \Lambda\land v''$,}
the maxmin asymptotic level~$\alpha$ test~$\delta_{q_v}=(\delta_{n,q_v})$
for~$H_v$ vs.~$K_v$ is given by
\begin{equation}
   \delta_{n,q_v} =  \lJc \biggl \lgroup \frac{1}{\sqrt{n}\,}\sum_{i=1}^n
   \Lambda^{(v)} (x_i) >  \Vert \Lambda^{(v)} \Vert
   \hspace{6\ei} u_{\alpha} + r_0 \hspace{6\ei}
   (v'' - v' \hspace{4\ei})   \biggr \rgroup
\end{equation}
and achieves maxmin asymptotic power
  $ \Phi \bigl(-u_{\alpha}+ \tau \hspace{6\ei}\Vert \Lambda^{(v)}
    \Vert  \,\bigr) $.
\item    \label{i.S.robT.v}
The test sequence~$\delta_{q_v}$ coincides with the robust asymptotic test
based on least favorable probability pairs
for~$U_v \bigl(P_{\theta_0};r_0/\!\sqrt{n}\,\bigr)$ vs.\/~{\rm
    $U_v \bigl(P_{\theta_0+\tau\!/\!\sqrt{n}\,};r_1/\!\sqrt{n}\,\bigr)$,}
hence maximizes the asymptotic minimum power over
  $U_v \bigl(P_{\theta_0+\tau\!/\!\sqrt{n}\,}; r_1/\!\sqrt{n}\,\bigr) $
subject to asymptotic maximum size~$\le\alpha$
over~$U_v \bigl(P_{\theta_0};r_0/\!\sqrt{n}\,\bigr) $.
\end{ABC}\end{Thm}
\newpage
\begin{Bew}\enskip
\begin{ABC}\item[{\rfi{\ref{i.S.maxminT.v}}}]
Also~$G_v$ is symmetric convex, so $ G_{10} =
  \tau \hspace{3\ei}\Lambda + r_1 \hspace{6\ei}G_v - r_0 \hspace{6\ei}G_v
  = \tau \hspace{3\ei}\Lambda - r \hspace{6\ei}G_v $,
and the minimum norm element~$q_{10}$ is supplied by
  $ q_0= r_0 \hspace{6\ei}\tilde{g} $,
  $ q_1= \tau \hspace{3\ei}\Lambda - r_1 \hspace{6\ei}\tilde{g} $,
where~$\tilde{g}\in G_v$ is the unique minimizer of
  $ \Vert \tau \hspace{3\ei}\Lambda-r \hspace{6\ei}g\Vert $
among all $g\in G_v$.
\par
The projection of~$\tau \hspace{3\ei}\Lambda$ on~$r \hspace{6\ei}G_h$
is determined by Theorem~\ref{t.rhov}, with~$\tau \hspace{3\ei}\Lambda$
in the place of~$\Lambda$. Then condition~(\ref{e.S.rad.v}) coincides
with condition~(\ref{e.rv}), and is equivalent to
  $r \hspace{6\ei}\tilde{g}\ne \tau \hspace{3\ei}\Lambda$,
that is, $G_{v,0}\cap G_{v,1}=\emptyset$. Thus~$\tilde{g}$
is of form~(\ref{e.S.clip.Lambda.v}),~(\ref{e.S.clip.c.v}).
\par
With $q_{10}=\tau \hspace{3\ei}\Lambda^{(v)}$
and  $ \langle q_{10}|q_0\rangle = \tau \hspace{6\ei} r_0
       \hspace{6\ei} (v'' - v' \hspace{4\ei})$,
Theorem~\ref{t.x} applies.
\item[{\rfi{\ref{i.S.robT.v}}}]
We invoke the results of HR~(1978), replacing $P_{-\tau_n}$ by~$P_0$
in~(2.8) there. This reduces~$2 \hspace{6\ei}\tau$ to~$\tau$ in that paper.
Then the radius condition~(2.6) of HR~(1978): $r/\tau<\Ew \Lambda_+$,
coincides with~(\ref{e.S.rad.v}). Moreover, the clipping equations~(3.9)
of HR~(1978) agree with~(\ref{e.S.clip.c.v}), and then the
function~$\Lambda^{(v)}$ equals the function~(3.10) of HR~(1978).
\par
Therefore, Theorems~3.4 and~4.1 of HR~(1978) tell us that
           $\delta_{q_v}$ maximizes the asymptotic minimum power
      over~$U_v (P_{\theta_0+s_n\tau};s_n \hspace{4\ei}r_1)$
subject to asymptotic maximum size~$\le\alpha$
      over~$U_v(P_{\theta_0};s_n \hspace{4\ei}r_0)$.\qed
\end{ABC}\end{Bew}
\begin{Rem}\rm  \label{r.S.nonu.lft} %
Under condition~(\ref{e.S.rad.v}), all least favorable pairs
  $g_v=(g_{v,0},g_{v,1})$ of tangents in~$G_{v,0}\times G_{v,1}$
are characterized by
\begin{equation}
      g_{v,0}= r_0 \hspace{6\ei}g_0\weg,\qquad
      g_{v,1}= \tau \hspace{3\ei}\Lambda -r_1 \hspace{6\ei}g_1
\end{equation}
where $g_0$ and~$g_1$ may be any elements of~$G_v$ whose positive
and negative parts make up those of~$\tilde{g}_v$ given
by~(\ref{e.S.clip.Lambda.v}) and~(\ref{e.S.clip.c.v}) such that
\begin{equation}
    r_0 \hspace{6\ei} g_0^+ + r_1 \hspace{6\ei} g_1^+ =
    \tau \hspace{6\ei}(\Lambda - v'' \hspace{6\ei})_+
\weg, \qquad
     r_0 \hspace{6\ei} g_0^- + r_1 \hspace{6\ei} g_1^- =
     \tau \hspace{6\ei}(v'- \Lambda)_+
\end{equation}
The least favorable tangent pair $q_v=(q_{v,0},q_{v,1})$, which results
from the special choice $ g_0 = g_1 = \tilde{g}_v $, is not the only one
in general. Other choices of $g_0$ and~$g_1$ may be based on suitable
partitions of the events $ \{ \Lambda> v'' \hspace{6\ei}\} $
                     and~$ \{ \Lambda< v'  \hspace{6\ei}\} $.
\par     For testing $ U_v \bigl(P_{\theta_0};r_0/\!\sqrt{n}\,\bigr) $
vs.~$ U_v \bigl(P_{\theta_0+\tau\!/\!\sqrt{n}\,};r_1/\!\sqrt{n}\,\bigr) $,
all least favorable pairs of probabilities have been characterized
by HR~(1977; Theorem~5.2).      
\qed\end{Rem}
\subsubsection{Maxmin Tests for Contamination Neighborhoods}
\begin{Thm}\sl \label{t.S.sad.c} Let
\begin{equation} \label{e.S.rad.c}
      r_0 < \hEw \bigl(\tau \hspace{3\ei}\Lambda-(r_1-r_0)\bigr)_+
\end{equation}
\begin{ABC}\item    \label{i.S.maxminT.c}
Then the least favorable tangent pair~$q_c=(q_{c,0},q_{c,1})$
     in~$G_{c,0}\times G_{c,1}$ is unique,
\begin{equation} \label{e.S.qc}  \mathsurround0em\arraycolsep0em
\hspace{-4em}
   q_{c,0} = \tau\hspace{6\ei}
             (\Lambda- c'' \hspace{6\ei})_+ -r_0 \weg, \qquad
   q_{c,1} = \left\{\begin{array}{l} \Ds  \tau \hspace{3\ei}\Lambda
                     + \tau \hspace{6\ei} (c'-\Lambda)_+ -r_1 \\
             \rule{0pt}{3ex} \Ds  \tau \hspace{6\ei}
                             (\Lambda\lor c' \hspace{6\ei}) -r_1
\end{array}\right.
\hspace{-3em}\end{equation}
with clipping constants
   $c'=c'(r_1/\tau)<z<c''(r_0/\tau)=c''$ determined by
\begin{equation}\label{e.S.clip.eta.c}
   \tau \hspace{6\ei}\hEw ( c' - \Lambda)_{+} = r_1 \weg,\qquad
   \tau \hspace{6\ei}\hEw ( \Lambda - c'' \hspace{6\ei})_{+} = r_0
\end{equation}
where $ z=(r_1-r_0)/\tau $.
Setting\/ {\rm $ \Lambda^{(c)} =  c'\lor \Lambda\land c'' - z $,}
the maxmin asymptotic level~$\alpha$ test~$\delta_{q_c}=(\delta_{n,q_c})$
for~$H_c$ vs.~$K_c$ is given by
\begin{equation}
   \delta_{n,q_c} =  \lJc \biggl \lgroup \frac{1}{\sqrt{n}\,}\sum_{i=1}^n
   \Lambda^{(c)} (x_i) >  \Vert \Lambda^{(c)} \Vert
   \hspace{6\ei} u_{\alpha} + r_0 \hspace{6\ei} (c'' - z) \biggr \rgroup
\end{equation}
and achieves maxmin asymptotic power
  $ \Phi \bigl(-u_{\alpha}+ \tau \hspace{6\ei}\Vert \Lambda^{(c)}\Vert
    \,\bigr) $.
\item    \label{i.S.robT.c}
The test sequence~$\delta_{q_c}$ coincides with the robust asymptotic test
based on least favorable probability pairs
for~$U_c \bigl(P_{\theta_0};r_0/\!\sqrt{n}\,\bigr)$ vs.\/~{\rm
    $U_c \bigl(P_{\theta_0+\tau\!/\!\sqrt{n}\,};r_1/\!\sqrt{n}\,\bigr)$,}
hence maximizes the asymptotic minimum power over
  $U_c \bigl(P_{\theta_0+\tau\!/\!\sqrt{n}\,}; r_1/\!\sqrt{n}\,\bigr) $
subject to asymptotic maximum size~$\le\alpha$
over~$U_c \bigl(P_{\theta_0};r_0/\!\sqrt{n}\,\bigr) $.
\end{ABC}\end{Thm}
\pagebreak[2]\begin{Bew}\enskip  
\begin{ABC}\item[{\rfi{\ref{i.S.maxminT.c}}}]       We can show that
  $ G_{10} = \tau \hspace{3\ei}\Lambda + r_1 \hspace{6\ei}G_c
                                       - r_0 \hspace{6\ei}G_c $
equals the closed set
\begin{equation}  \label{e.s.Gc-Gc=Gv}
    \tau \hspace{3\ei}\Lambda -(r_1-r_0)
    + \bigl\{\,g\in L_2 \bigm| \Ew g = r_1-r_0\weg,
    \:\Ew g^+\le r_1\weg,\:\Ew g^-\le r_0 \,\bigr\}
\end{equation}
As $\Ew (\tau \Lambda-c)_+ = \Ew (c- \tau \Lambda)_+ - c $,
radius condition~(\ref{e.S.rad.c}) is equivalent to
\begin{equation} \label{e.S.rad.c1}
  r_1 < \hEw \bigl((r_1-r_0)-\tau \hspace{3\ei}\Lambda\bigr)_+
\end{equation}
If (\ref{e.S.rad.c}) and~(\ref{e.S.rad.c1}) are violated,
the zero function is in~$G_{10}$ as
$$ 0 = \tau \hspace{3\ei}\Lambda -(r_1-r_0)
     + \bigl((r_1-r_0)-\tau \hspace{3\ei}\Lambda\bigr)_+
     - \bigl(\tau \hspace{3\ei}\Lambda-(r_1-r_0)\bigr)_+     $$
Under conditions (\ref{e.S.rad.c}) and~(\ref{e.S.rad.c1}),
equivalently $ c'<z=(r_1-r_0)/\tau<c''$ for the solutions $c'$
and~$c''$ to~(\ref{e.S.clip.eta.c}), the function
  $ q_{10}=q_{c,1}-q_{c,0} $ is nonzero,
\begin{equation} \label{e.s.q10.c}
\mathsurround0em\arraycolsep0em \begin{array}{r@{{}={}}l}
\Ds q_{10} & \Ds     \tau \hspace{3\ei}\Lambda -(r_1-r_0) +
                     \tau \hspace{6\ei} (c'-\Lambda)_+ -
   \tau\hspace{6\ei} (\Lambda- c'' \hspace{6\ei})_+ \\
\rule{0pt}{3ex} & \Ds
    \tau \hspace{6\ei}(c'\lor \Lambda\land c'' \hspace{6\ei})
    -(r_1-r_0) = \tau \hspace{3\ei} \Lambda^{(c)}
\end{array}\end{equation}
and, by Lemma~\ref{l.c.approx}, the minimum norm element of~$G_{10}$.
In fact, for all $g_0\in G_c$,
\begin{equation} \label{e.s.gLeta.Dg0}
\mathsurround0em\arraycolsep0em \begin{array}{rcl}
\Ds  \langle \Lambda^{(c)}| r_0 \hspace{6\ei}g_0- q_{c,0}\rangle
     & {}={} & \Ds \bigl \langle c'\lor \Lambda\land c'' \big|
        r_0 \hspace{6\ei}(1+g_0) - \tau\hspace{6\ei}
        (\Lambda- c'' \hspace{6\ei})_+ \bigr\rangle    \\
\rule{0pt}{3ex} & {}\le{} &  \Ds
        c'' \hspace{3\ei}r_0 \hspace{3\ei}\Ew(1+g_0) -
        c'' \hspace{3\ei}\tau \hspace{3\ei}
        \Ew (\Lambda- c'' \hspace{6\ei})_+ = 0
\end{array}\end{equation}
as $ c'\lor \Lambda\land c'' \le c'' $ and $ 1+g_0 \ge0 $,
and by~(\ref{e.S.clip.eta.c}).
Likewise, for all $g_1\in G_c$,
\begin{equation} \label{e.s.gLeta.Dg1}
\mathsurround0em\arraycolsep0em \begin{array}{rcl}
\Ds \langle \Lambda^{(c)}| q_{c,1}
      - \tau \hspace{3\ei}\Lambda - r_1 \hspace{6\ei}g_1 \rangle
    & {}={} & \Ds   \bigl \langle c'\lor \Lambda\land c'' \big|
        \tau\hspace{6\ei} (c'-\Lambda)_+
        - r_1 \hspace{6\ei}(1+g_1) \bigr\rangle    \\
\rule{0pt}{3ex} &{}\le{}&\Ds
        c' \hspace{3\ei}\tau \hspace{3\ei} \Ew (c' -\Lambda)_+
        - c' \hspace{3\ei}r_1 \hspace{3\ei}\Ew(1+g_1) = 0
\end{array}\samepage\end{equation}
With $ \langle \Lambda^{(c)}| q_{c,0}\rangle =
        r_0 \hspace{6\ei}(c''-z)$, Theorem~\ref{t.x} applies.
\par
Now let 
$ (r_0 \hspace{6\ei}g_0, \tau\hspace{3\ei} \Lambda + r_1 \hspace{6\ei}g_1) $
be any least favorable tangent pair, that is, with elements $g_0,g_1\in G_c$
such that $ \tau \hspace{3\ei}\Lambda + r_1 \hspace{6\ei}g_1
                                     - r_0 \hspace{6\ei}g_0 = q_{10} $.
Then, in view of~(\ref{e.s.q10.c}),
\begin{equation}
       r_1 \hspace{6\ei}(1+g_1) - r_0 \hspace{6\ei}(1+g_0) =
       \tau \hspace{6\ei} (c'-\Lambda)_+ -
       \tau\hspace{6\ei} (\Lambda- c'' \hspace{6\ei})_+
\end{equation}
The RHS, since $c'<c''$, is a decomposition into positive and
negative parts. As also $1+g_1\ge0$ and $1+g_0\ge0$ this implies that
\begin{equation}
       \tau \hspace{6\ei} (c'-\Lambda)_+ \le r_1 \hspace{6\ei}(1+g_1)
\weg,\qquad   \tau\hspace{6\ei} (\Lambda- c'' \hspace{6\ei})_+ \le
              r_0 \hspace{6\ei}(1+g_0)
\end{equation}
But by~(\ref{e.S.clip.eta.c}), the functions compared have the same
expectations. Hence strict inequalities cannot hold. It follows that
\begin{equation}
        r_0 \hspace{6\ei}g_0 =
        \tau\hspace{6\ei} (\Lambda- c'' \hspace{6\ei})_+ -r_0
\weg, \qquad
        r_1 \hspace{6\ei}g_1 =
        \tau \hspace{6\ei} (c'-\Lambda)_+ -r_1
\end{equation}
which proves uniqueness of the least favorable tangent
        pair~$q_c=(q_{c,0},q_{c,1})$.
\item[{\rfi{\ref{i.S.robT.c}}}]
The substitution of $P_{-\tau_n}$ by~$P_0$ in HR~(1978)
reduces~$2 \hspace{6\ei}\tau$ to~$\tau$ there. Then the radius
condition~(2.6) of HR~(1978) is~(\ref{e.S.rad.c}). Moreover, the clipping
equations~(3.9) of HR~(1978) agree with~(\ref{e.S.clip.eta.c}), and the
present function~$\Lambda^{(c)}$ 
equals the function defined by~(3.10) in HR~(1978).
\par
Therefore, Theorems~3.4 and~4.1 of HR~(1978) tell us that~$\delta_{q_c}$
maximizes the asymptotic minimum power
over~$U_c (P_{\theta_0+s_n\tau};s_n \hspace{4\ei}r_1)$ subject to
asymptotic maximum size~$\le\alpha$
      over~$U_c(P_{\theta_0};s_n \hspace{4\ei}r_0)$.\qed
\end{ABC}\end{Bew}
\begin{Rem}\rm   The radius condition~(\ref{e.S.rad.c}),
being equivalent to
    $ \tau \hspace{6\ei}c''> r_1-r_0 $
for~$c''$ satisfying~(\ref{e.S.clip.eta.c}),
is stronger than
    $ \tau \hspace{6\ei}c''> -r_0 $.
In turn,  $ \tau \hspace{6\ei}c''> -r_0 $
for~$c''$ satisfying~(\ref{e.S.clip.eta.c}),
can be shown to be equivalent to
    $ r_0 < - \tau \hspace{6\ei}\winf_P \Lambda $.
\par     Under this radius condition~(\ref{e.rc}):
         $ r_0 < - \tau \hspace{6\ei}\winf_P \Lambda $, Theorem~\ref{t.rhoc}
(with $\tau \hspace{3\ei}\Lambda$ in the place of~$\Lambda$) yields the
element~$\tilde{g}_0$ of~$G_c$ minimizing
    $ \Vert \tau \hspace{3\ei}\Lambda - r_0 \hspace{6\ei}g\Vert $
among all $g\in G_c$:
\begin{equation}
  r_0 \hspace{6\ei}\tilde{g}_0 =
  \tau \hspace{3\ei} \Lambda - (\tau \hspace{3\ei} \Lambda + r_0)\land u =
  \tau \hspace{6\ei} (\Lambda- c'' \hspace{6\ei})_+ - r_0    
\end{equation}
with $u$ and~$\tau \hspace{6\ei}c''= u-r_0 $ determined by
   $\Ew \tilde{g}_0=0$. Thus, $ q_{c,0}= r_0 \hspace{6\ei} \tilde{g}_0 $.
\par     Likewise, the radius condition~(\ref{e.S.rad.c1}),
being equivalent to
    $ \tau \hspace{6\ei}c'< r_1-r_0 $
for~$c'$ satisfying~(\ref{e.S.clip.eta.c}),
implies that $ \tau \hspace{6\ei}c' < r_1 $, equivalently
    $ r_1 < \tau \hspace{6\ei}\wsup_P \Lambda $.
\par     Under this radius condition~(\ref{e.rc}):
         $ r_1 < \tau \hspace{6\ei}\wsup_P \Lambda $, Theorem~\ref{t.rhoc}
(with $-\tau \hspace{3\ei}\Lambda$ in the place of~$\Lambda$) yields
the element~$\tilde{g}_1$ of~$G_c$ minimizing
    $ \Vert \tau \hspace{3\ei}\Lambda + r_1 \hspace{6\ei}g\Vert $
among all $g\in G_c$. And then it may again be verified that
    $ q_{c,1}= \tau \hspace{3\ei}\Lambda + r_1 \hspace{6\ei} \tilde{g}_1 $.
\par  Therefore, according to Lemma~\ref{l.c.approx}, it follows that,
for all $g_0,g_1\in G_c$,
\begin{equation}
    \langle \tau \hspace{3\ei}\Lambda -r_0 \hspace{6\ei}\tilde{g}_0|
    r_0 \hspace{6\ei}g_0 - r_0 \hspace{6\ei}\tilde{g}_0 \rangle \le 0
\weg, \qquad
    \langle \tau \hspace{3\ei}\Lambda +r_1 \hspace{6\ei}\tilde{g}_1|
    r_1 \hspace{6\ei}\tilde{g}_1 - r_1 \hspace{6\ei}g_1 \rangle \le 0
\end{equation}
But the bounds (\ref{e.s.gLeta.Dg0}) and~(\ref{e.s.gLeta.Dg1}) established
in the preceeding proof tell us that this remains true for
  $ \tau \hspace{3\ei} \Lambda^{(c)} 
    = \tau \hspace{3\ei}\Lambda + r_1 \hspace{6\ei}\tilde{g}_1
      - r_0 \hspace{6\ei}\tilde{g}_0 $ in the place of
   $ \tau \hspace{3\ei}\Lambda -r_0 \hspace{6\ei}\tilde{g}_0 $, respectively
of $ \tau \hspace{3\ei}\Lambda +r_1 \hspace{6\ei}\tilde{g}_1 $.
This is remarkable since the two additional terms are always nonnegative,
\begingroup \mathsurround0em\arraycolsep0em \begin{eqnarray}
\hspace{-1.5em}  \langle r_1 \hspace{6\ei}\tilde{g}_1 |
  r_0 \hspace{6\ei}g_0 - r_0 \hspace{6\ei}\tilde{g}_0 \rangle
  &{}={}&
  \bigl \langle \tau \hspace{6\ei} (c'-\Lambda)_+ -r_1\big|
    r_0 \hspace{6\ei}g_0+r_0- \tau \hspace{6\ei}
                 (\Lambda- c'' \hspace{6\ei})_+ \bigr \rangle
  \nonumber \\ &{}={}&
  \tau \hspace{6\ei}r_0 \hspace{6\ei}
      \langle (c'-\Lambda)_+| 1+g_0\rangle \ge0 \\
\hspace{-1.5em}  \langle r_0 \hspace{6\ei}\tilde{g}_0 |
  r_1 \hspace{6\ei}g_1 - r_1 \hspace{6\ei}\tilde{g}_1 \rangle
  &{}={}&
  \bigl \langle \tau \hspace{6\ei} (\Lambda- c'' \hspace{6\ei})_+
  -r_0  \big| r_1 \hspace{6\ei}g_1+r_1- \tau \hspace{6\ei}
                 (c'-\Lambda)_+ \bigr \rangle
  \nonumber \\ &{}={}&
  \tau \hspace{6\ei}r_1 \hspace{6\ei}
    \langle (\Lambda- c'' \hspace{6\ei})_+| 1+g_1\rangle \ge0
\end{eqnarray}\endgroup
where use has been made of $ c'< c''$, which is garanteed by the
stronger radius condition~(\ref{e.S.rad.c}),~(\ref{e.S.rad.c1}).
\qed\end{Rem}
\begin{Rem}\rm    For testing
    $ U_c \bigl(P_{\theta_0};r_0/\!\sqrt{n}\,\bigr) $
vs.~$ U_c \bigl(P_{\theta_0+\tau\!/\!\sqrt{n}\,};r_1/\!\sqrt{n}\,\bigr) $,
all least favorable pairs of probabilities have been characterized
in terms of their densities by HR~(1977; Theorem~5.2).
The uniqueness of the least favorable tangent pair~$q_c=(q_{c,0},q_{c,1})$
gives rise to the conjecture that, contrary to the total variation case,
\begin{equation}
   \lim_{n\to \infty}\sqrt{n}\:d_h (Q''_{n,j}, Q'_{n,j})=0 \weg,
   \qquad j=0,1 \hspace{-2.25em}
\end{equation}  
if  $(Q'_{n,0},Q'_{n,1})$ and $(Q''_{n,0},Q''_{n,1})$
are any two, possibly different, least favorable probability pairs
for $ U_c \bigl(P_{\theta_0};r_0/\!\sqrt{n}\,\bigr) $
vs.~$ U_c \bigl(P_{\theta_0+\tau\!/\!\sqrt{n}\,};r_1/\!\sqrt{n}\,\bigr) $.
\qed\end{Rem}
\begin{Rem}\rm     For shrinking contamination neighborhoods of a one
parameter family involving a finite dimensional nuisance parameter,
the robust asymptotic tests based on least favorable pairs were investigated
by Wang~(1981).    It would be interesting to derive his maxmin asymptotic
test by projection.
\qed\end{Rem}
\paragraph{Acknowledgement}
A first version of the paper was written during a visit to
Sonderforschungsbereich~373
      (``Quantifikation und Simulation \"okonomischer Prozesse'')
      at Humboldt--Universit\"at zu Berlin.
Financial support by Deutsche Forschungsgemeinschaft is acknowledged.
\par     I thank Wolfgang H\"ardle for his hospitality.
I am indebted to Peter Ruckdeschel for his assistance and
the figures. 
\sz 

\nz 
\vfil
\settowidth{\breit}{\sc  e-mail:
                    \rm helmut.rieder@uni-bayreuth.de}
\hfill
\parbox{\breit}{\sc 
                    Department of Mathematics\\
                    University of Bayreuth, NW~II\\
                    D-95440 Bayreuth, Germany\\   e-mail:
                    {\rm helmut.rieder@uni-bayreuth.de}}

\begin{thebibliography}{99} \frenchspacing
\bibitem{} Beran, R.J. (1974):\hskip.5em Asymptotically efficient and
       adaptive rank estimates in location models.
       {\sl Ann. Statist.\/} {\bf 2} 63--74.
\bibitem{} \mbox{Bickel, P.J.~(1981)}:\hskip.5em
   Quelques aspects de la statistique robuste.
   In {\sl Ecole d'Et\'e de Probabilit\'es de Saint Flour IX 1979\/}
   (P.L.~Hennequin, ed.), 1--72.
   Lecture Notes in Mathematics~\#876. Springer, Berlin.
\bibitem{} \mbox{Bickel, P.J.~(1982)}:\hskip.5em On adaptive estimation.
   {\sl Ann. Statist.\/} {\bf 10} 647--671.
\bibitem{} \mbox{Bickel, P.J.}, Klaassen, C.A.J., Ritov, Y.,
            and Wellner, J.A.~(1993):\hskip.5em
  {\sl Efficient and Adaptive Estimation for Semiparametric Models\/}.
  Springer, New York.
\bibitem{} Birg\'e, L. (1980): Approximation dans les espaces m\'etrique
   et th\'eorie de l'estimation. Ph.D.~Thesis, University of Paris.
\bibitem{} Hampel, F.R., Ronchetti, E.M., Rousseeuw, P.J. and
         Stahel, W.A.~(1986):\hskip.5em
         {\sl \mbox{Robust} Statistics---The Approach Based on
         Influence Functions\/}. Wiley, New York.
\bibitem{} \mbox{Huber, P.J.~(1964)}:\hskip.5em
     Robust estimation of a location parameter.
     {\sl Ann. Math. Statist.\/} {\bf 35} 73--101.
\bibitem{} \mbox{\rm Huber,~P.J.~(1981)}:\hskip.5em
           {\sl Robust Statistics\/}. Wiley, New York.
\bibitem{} \mbox{Huber, P.J.~(1996)}:\hskip.5em
       {\sl Robust Statistical Procedures\/} (2nd ed.). CBMS-NSF
       Regional Conference Series in Applied Mathematics \#~68.
       Soc. Industr. Appl. Math., Philadelphia, Pennsylvania.
\bibitem{} \mbox{Huber,~P.J.} and Strassen,~V.(1973):\hskip.5em
        Minimax tests and the Neyman--Pearson lemma for capacities.
        {\sl Ann.\ Statist.\/}~{\bf1} \mbox{251--263}.
\bibitem{} Huber--Carol, C.~(1970):\hskip.5em
         {\sl \'Etude asymptotique de tests robustes\/}.
         Th\`ese de \mbox{Doctorat}, ETH~Z\"urich.
\bibitem{} Kakiuchi, I. and Kimura, M.~(2000):
           Robust rank tests for $k$-sample approximate equality
           in the presence of gross errors.
           {\sl J. Statist. Plann. Inf.\/}~{\bf 82}, to appear.
\bibitem{} Klaassen, C.A.J. (1987):\hskip.5em
        Consistent estimation of the influence function of locally
        asymptotically linear estimators.
        {\sl Ann.\ Statist.\/} {\bf15} \mbox{1548--1562}.
\bibitem{} \mbox{Pfanzagl, J.} and Wefelmeyer, W.~(1982):\hskip.5em
  {\sl Contributions to a General Asymptotic Statistical Theory\/}.
  Lecture Notes in Statistics~\#13. Springer, Berlin.
\bibitem{}   Quang, P.X.~(1985):\hskip.5em
             Robust sequential testing.
             {\sl Ann.\ Statist.\/}~{\bf13} \mbox{638--649}.
\bibitem{} \mbox{\rm Rieder,~H.~(1977)}:\hskip.5em
        Least favorable pairs for special capacities.
        {\sl Ann.\ Statist.\/} {\bf5} \mbox{909--921}.
\bibitem{} \mbox{\rm Rieder,~H.~(1978)}:\hskip.5em
        A robust asymptotic testing model.
        {\sl Ann.\ Statist.\/}~{\bf6} \mbox{1080--1094}.
\bibitem{} \mbox{Rieder,~H.~(1980)}:\hskip.5em
        Estimates derived from robust tests.
        {\sl Ann.\ Statist.\/}~{\bf8} \mbox{106--115}.
\bibitem{} \mbox{Rieder,~H.~(1981\,a)}:\hskip.5em
       Robustness of one- and two-sample rank tests against gross errors.
       {\sl Ann. Statist.\/}~{\bf 9} 245--265.
\bibitem{} \mbox{Rieder,~H.~(1981\,b)}:\hskip.5em
    On local asymptotic minimaxity and admissibility in robust estimation.
    {\sl Ann.\ Statist.\/}~{\bf 9} \mbox{266--277.}
\bibitem{} \mbox{Rieder,~H.~(1994)}:\hskip.5em
           {\sl Robust Asymptotic Statistics\/}. Springer, New York.
\bibitem{} \mbox{Rieder,~H.~(2000)}:\hskip.5em
   One-sided confidence about functionals over tangent cones.
   Submitted for publication.
\bibitem{} \mbox{Shen, L.Z.~(1994)}:\hskip.5em
   Optimal robust estimates for semiparametric symmetric location models.
   {\sl Statistics\&Decisions\/}~{\bf 12} 113--124.
\bibitem{} \mbox{Shen, L.Z.~(1995)}:\hskip.5em
   On optimal B-robust influence functions in semiparametric models.
   {\sl Ann. Statist.\/}~{\bf 23} 968--989.
\bibitem{}  Schick, A. (1986):\hskip.5em
   On asymptotically efficient estimation in semiparametric models.
   {\sl Ann. Statist.\/}~{\bf 14} 1139--1151.
\bibitem{} Stone, C.~(1975):\hskip.5em
   Adaptive maximum likelihood estimation for a location parameter.
   {\sl Ann. Statist.\/} {\bf 3} 267--284.
\bibitem{} van der Vaart, A.W.~(1998):\hskip.5em
   {\sl Asymptotic Statistics\/}. CUP, Cambridge.
\bibitem{} Wang, P.C.C.~(1981):\hskip.5em
   Robust asymptotic tests of statistical hypotheses involving nuisance
   parameters.   {\sl Ann.\ Statist.\/}~{\bf9} \mbox{1096--1106}.
\nonfrenchspacing \end{thebibliography}
\end{document}